\documentclass[11pt,reqno]{amsart}

\usepackage{amsfonts,amsmath,amsthm}
\usepackage{amssymb,epsfig}
\usepackage{enumerate} 

\usepackage{algorithm}
\usepackage{algpseudocode}
\usepackage{hyperref}
\usepackage{color} %color
\definecolor{vert}{rgb}{0,0.6,0}

\usepackage[text={425pt,650pt},centering]{geometry}
\usepackage{graphicx}
\usepackage{epsfig}
\usepackage{tikz,esvect}
\usetikzlibrary{3d,calc,shapes}
\usepackage{pgfplots}
\pgfplotsset{compat=1.18}
\usepackage{booktabs}

\usepackage{comment}

\usepackage{caption}
\usepackage{color} %color
\definecolor{vert}{rgb}{0,0.6,0}

\numberwithin{figure}{section}

\theoremstyle{plain}
\newtheorem{thm}{Theorem}[section]
\newtheorem{ass}{Assumption}
\newtheorem{prob}{Problem}
\newtheorem{defn}{Definition}

\newtheorem{lem}[thm]{Lemma}
\newtheorem{cor}[thm]{Corollary}
\newtheorem{prop}[thm]{Proposition}
\theoremstyle{remark}
\newtheorem{rem}{\bf{Remark}}
\numberwithin{equation}{section}
\theoremstyle{definition}

\usepackage{etoolbox}
\newtheorem{hjbsetting}{}

\AtEndEnvironment{hjbsetting}{\hfill\(\Diamond\)\par}

\usepackage{ulem}

\newcommand{\N}{\mathbb{N}}

\newcommand{\R}{\mathbb{R}}

\newcommand{\T}{\mathbb{T}}

\newcommand{\cS}{\mathcal{S}}

\newcommand{\cE}{\mathcal{E}}

\newcommand{\BUC}{{\rm BUC\,}}

\newcommand{\ep}{\varepsilon}

\newcommand{\Lop}{\mathcal{L}}
\newcommand{\norm}[1]{\left\|#1\right\|}

\newcommand{\Dp}{D_p}
\newcommand{\Dpp}{D^2_{pp}}
\newcommand{\argmin}{\operatorname{argmin}}
\newcommand{\argmax}{\operatorname{argmax}}

\begin{document}

\title[Newton's Method and Policy Iterations for Viscous HJ Equations]{Newton's Method and Policy Iterations for Viscous Hamilton--Jacobi Equations}
\author[D. A. Gomes, X. Guo, Y. Kim, H. V. Tran, Y. P. Zhang]{Diogo A. Gomes, Xiaoqin Guo, Yeoneung Kim, Hung V. Tran, Yuming P. Zhang}

\thanks{
D. A. Gomes is supported by King Abdullah University of Science and Technology (KAUST) baseline funds.
X. Guo is supported by the Simons Foundation through a Collaboration Grant for Mathematicians \#852943. 
Y. Kim is supported by the National Research Foundation of Korea (NRF) grants funded by the Korea government (MSIT) (RS-2023-00219980, RS-2026-25490291). H. V. Tran is partially supported by NSF grant DMS-2348305. 
Y. P. Zhang is partially supported by NSF CAREER grant DMS-2440215, Simons Foundation Travel Support MPS-TSM-00007305.
}

\address[D. A. Gomes]
{King Abdullah University of Science and Technology (KAUST), CEMSE Division, Thuwal 23955-6900, Saudi Arabia}
\email{diogo.gomes@kaust.edu.sa}

\address[X. Guo]
{
Department of Mathematical Sciences, 
University of Cincinnati,
2815 Commons Way, Cincinnati OH 45221-0025}
\email{guoxq@ucmail.uc.edu}

\address[Y. Kim]
{Department of Industrial Engineering, 
Yonsei University, 50 Yonsei-ro, Seodaemun-gu, Seoul 03722, Republic of Korea}
\email{yeoneung@yonsei.ac.kr}

\address[H. V. Tran]
{
Department of Mathematics, 
University of Wisconsin-Madison, Van Vleck Hall, 480 Lincoln Drive, Madison, WI 53706}
\email{hung@math.wisc.edu}

\address[Y. P. Zhang]
{
Department of Mathematics and Statistics, Auburn University, STEM \& AG Building, 305 W. Samford Avenue, Auburn, AL 36849}
\email{yzhangpaul@auburn.edu}

\date{}


% 35K55: Nonlinear parabolic equations.
% 65J15: Numerical solutions to equations with nonlinear operators.
% 65M15: Error bounds for initial value and initial-boundary value
%        problems involving PDEs.

\begin{abstract}
We study Newton's method for the Cauchy problem for viscous Hamilton--Jacobi equations with Hamiltonians that need not be convex. 
Each Newton step solves a linear parabolic equation obtained by linearizing the equation at the current iterate. 
We prove that the iterates converge geometrically in
\(L^\infty\) from any initial guess satisfying the initial condition, and quadratically in the parabolic H\"older space \(C^{2+\alpha,1+\alpha/2}\) from initial guesses close to the solution. 
When the Hamiltonian admits a control representation,
the Newton iterates coincide with those of classical policy iteration, which therefore inherits both convergence rates. We extend the local quadratic convergence to uniformly parabolic fully nonlinear equations. 
For small viscosity, we introduce a Nash--Moser variant that smooths the iterate and the residual in space, and we give a convergence criterion with constants independent of the viscosity. 
Finally, we show how computable residuals certify the error of inexact Newton corrections computed with physics-informed neural networks.
\end{abstract}

\maketitle

\section{Introduction}

% Background and importance of the field.
Viscous Hamilton--Jacobi equations arise in stochastic optimal control, large deviations,
differential games, and mean-field games; see, e.g.,
\cite{Fleming2006, MTL17, lasryMeanFieldGames2007}. Except in special cases, their solutions
must be computed numerically, and iterative schemes that reduce the
nonlinear problem to a sequence of linear ones are a central tool for doing
so. The best known of these is policy iteration, introduced by Howard
\cite{Howard1960DynamicProgramming}: it alternates between evaluating a fixed policy, which
requires solving a linear equation, and improving that policy pointwise.
For Markov decision processes, policy iteration is equivalent to Newton's
method applied to the Bellman equation \cite{putermanConvergencePolicyIteration1979}. This
equivalence explains its fast convergence and suggests extending
Newton-type methods beyond the convex setting of control theory.

% Problem, motivation, and setting.

We study approximation schemes for the following viscous Hamilton--Jacobi
Cauchy problem. Fix \(\alpha\in(0,1)\) and \(T>0\), and set
\[
    Q_T:=\R^n\times(0,T),
    \qquad
    \overline Q_T:=\R^n\times[0,T],
    \qquad
    X_T:=C_b^{2+\alpha,1+\alpha/2}(\overline Q_T).
\]
\begin{hjbsetting}
\label{set:HJB}
With the notation above, let \(H:\R^n\times\R^n\to\R\) and let
\[
    g\in C_b^{2+\alpha}(\R^n).
\]
We consider the Cauchy problem
\begin{equation}\label{eq:HJB}
\begin{cases}
    u_t+H(x,Du)=\Delta u
    \qquad & \text{in } Q_T,\\
    u(x,0)=g(x)
    \qquad & \text{on } \R^n.
\end{cases}
\end{equation}
Under the structural and regularity assumptions stated in
Section~\ref{sec:assumptions}, \(u\in X_T\) denotes the unique classical
solution of \eqref{eq:HJB}.
\end{hjbsetting}

The well-posedness theory for \ref{set:HJB} is not the main subject
of this paper.  Our goal is to approximate the solution \(u\) by Newton-type,
policy-iteration-type, and inexact computational schemes, and to prove
quantitative convergence and certification results for these approximations.
This leads to the following approximation problem.

\begin{prob}[Approximation and certification]
\label{prob:approx-HJB}
Consider the setting of \ref{set:HJB}, and let \(u\in X_T\) be the
corresponding solution. Construct a sequence
\(\{u^k\}_{k\ge0}\subset X_T\) with
\[
    u^k(x,0)=g(x)
    \qquad\text{for all } k\ge0,
\]
such that \(u^k\to u\) in \(L^\infty(\overline Q_T)\), and identify
conditions under which the convergence also holds in \(X_T\), with
quantitative rates.
In addition, develop numerical versions of the
iteration for which the approximation error can be certified by computable
residuals.
\end{prob}

% Previous work and gap in the literature.

A key feature of our analysis is that \(H(x,p)\) need not be convex in $p$.
%%%GGG:it may be better to metion p here, because D_pH shows up below.
Policy-iteration methods for nonconvex viscous Hamilton--Jacobi equations
were studied by Guo--Tran--Zhang \cite{Guo-Tran-Zhang}. Our approach instead linearizes the equation itself, without using a min--max
representation or its differential-game interpretation; in return, we require
suitable smoothness of \(H\).

% More concrete description of the problem addressed here.
To address Problem~\ref{prob:approx-HJB}, we apply Newton's method directly
to the nonlinear operator
\[
    v\longmapsto v_t+H(x,Dv)-\Delta v.
\]
At each step, we linearize this operator at \(u^k\) and solve the resulting
linear uniformly parabolic equation for $u^{k+1}$. The increment \(h^k:=u^{k+1}-u^k\), called the Newton correction, satisfies
\(h^k(\cdot,0)=0\).
The precise definition and well-posedness of the iteration are given in
Subsection~\ref{ss2.3}, while its relation to classical policy iteration is
established in Section~\ref{sec:policy-iteration}.

% Main results and how they improve on previous work.

Our first two main results give complementary convergence estimates for the
Newton iteration. Theorem~\ref{thm:main1} shows that, from every admissible
initial guess $u^0$,
\[
    \|u^k-u\|_{L^\infty(\overline Q_T)}
    \leq
    C\|D(u^0-u)\|_{L^\infty(\overline Q_T)}\,2^{-k},
    \qquad k\geq1,
\]
where \(C\) depends on \(T\) and on
\(\|D_pH\|_{L^\infty(\R^n\times\R^n)}\) for the modified Hamiltonian of
Remark~\ref{truncatedham}. 
The geometric convergence rate obtained here is related to earlier quantitative convergence results for policy iteration in the convex \cite{kerimkulov2020exponential, tang2025policy, kim2025neural, bensoussan2026operatorsplittingpolicyiteration,yang2025gaussian,yao2026model} and nonconvex \cite{ito2021neural, Guo-Tran-Zhang, yang2025solving} settings. Policy iteration has also been studied for exploratory Hamilton--Jacobi--Bellman equations \cite{huang2025convergence,ma2026convergence,tran2025policy}. Related developments include applications to periodic Hamilton--Jacobi--Bellman--Isaacs equations and numerical homogenization \cite{Kawecki-Sprekeler}, as well as neural implementations for Hamilton--Jacobi--Bellman equations \cite{lee2025hamilton,kim2025neural}. For related work on Newton iterations for inverse problems involving Hamilton--Jacobi--Bellman equations, see \cite{le2022carleman, abhishek2024carleman}.

Moreover, Theorem~\ref{thm:main2} shows that, if
\[
    \delta:=\|u^0-u\|_{X_T}\leq\delta_0,
\]
then
\[
    \|u^k-u\|_{X_T}
    \leq
    C_*^{-1}(C_*\delta)^{2^k},
    \qquad k\geq0.
\]
Thus, the iteration converges globally and geometrically in \(L^\infty\),
while sufficiently accurate initial guesses yield quadratic convergence in
the stronger parabolic Schauder norm \(X_T\). 
Smears and S\"uli \cite{Smear-Suli} analyzed a locally superlinearly convergent semismooth Newton method for finite element discretizations of Hamilton--Jacobi--Bellman equations with Cordes coefficients. 
In the context of Markov decision processes, a quadratic convergence rate of policy iteration has also been proven under certain conditions \cite{putermanConvergencePolicyIteration1979,santos2004convergence}.
For the entropy-regularized control problem, Ma--Wang--Zhang \cite{ma2026convergence} obtained quadratic convergence of policy iteration via the Feynman--Kac formula and a probabilistic approach.

The same local quadratic mechanism extends beyond the semilinear equation
\eqref{eq:HJB}: Theorem~\ref{thm:main3} establishes the corresponding result
for uniformly parabolic fully nonlinear equations. For the equation with viscosity \(\ep\Delta u\), \(0<\ep\le1\), a different
argument is needed, since the parabolic Schauder estimates are not uniform
as \(\ep\to0\). We therefore smooth both the
current iterate and the residual in space, and
Theorem~\ref{thm:criterion-new} gives a conditional Nash--Moser-type
criterion yielding residual decay and convergence of the resulting
iteration.

% Key ideas and techniques.

The proof of the global \(L^\infty\) estimate does not rely on smallness
of the initial error. Instead, we apply a comparison argument to
\(\tfrac12|u^k-u|^2\) and introduce auxiliary heat potentials, driven by
the gradient errors, that contract from one iteration to the next.
 The local quadratic estimates follow
from the integral Taylor formula for the Newton remainder combined with
global parabolic Schauder estimates, 
whose constants are uniform while the iterates remain in a fixed
\(X_T\)-ball around the solution.
In the
small-viscosity scheme, the spatial smoothing uses kernels whose moments
vanish up to high order, and we combine transport estimates
uniform in the viscosity with heat-kernel H\"older estimates whose
dependence on the viscosity is explicit. These ingredients balance in a
Nash--Moser scheme with doubly exponentially growing smoothing parameters.

% Contribution in computational aspect

From a computational viewpoint, the Newton formulation also leads naturally to an a posteriori error estimate when the linearized equations are solved approximately. After each approximate linearized step, the error of the updated iterate is controlled by the residual of that step together with the quadratic Newton remainder. A useful feature of this estimate is that residual errors from earlier iterations do not accumulate explicitly: their effect is already contained in the current iterate and is reassessed through its current residual (Theorem \ref{thm:pinn-aposteriori}). 
The numerical experiments illustrate this behavior for different training schedules and nonconvex high-dimensional problems, and show that the residual indicator remains informative when training stagnates or fails to resolve the solution.
These experiments further demonstrate the potential of neural-network-based Newton iterations for high-dimensional nonconvex problems, combining efficient linearized computations with empirical residual-based error monitoring in dimensions up to \(100\). The code and data for these experiments are publicly available.\footnote{\url{https://github.com/yeoneung/newton_policy}}

% Organization of the paper.
The paper is organized as follows. Section~\ref{sec:na} introduces the
notation and assumptions, and defines the Newton iteration.
Section~\ref{sec:newton-convergence} proves the global \(L^\infty\)
convergence (Theorem~\ref{thm:main1}) and the local quadratic convergence
(Theorem~\ref{thm:main2}). Section~\ref{sec:policy-iteration} relates the
Newton iteration to classical policy iteration.
Section~\ref{sec:fully-nonlinear} extends the local quadratic convergence
to fully nonlinear equations (Theorem~\ref{thm:main3}).
Section~\ref{sec:small-viscosity} develops the Nash--Moser iteration for
small viscosity (Theorem~\ref{thm:criterion-new}). Sections~\ref{sec:pinn-certificate} and~\ref{sec:numerics} present
numerical experiments and the residual-based certification of inexact Newton corrections.

\section{Setting and the Newton Iteration}
\label{sec:na}

In this section, we introduce the spatial and parabolic H\"older norms
used throughout the paper and state the assumptions on \(H\). These
assumptions yield a global gradient bound for the solution of
\eqref{eq:HJB}, so \(H\) may be modified for large \(|p|\) without
changing the equation satisfied by \(u\). We then define the Newton
iteration, show that each step is well posed, and identify it as the
Newton linearization of the PDE residual. For convex Hamiltonians, this
linearization is the policy-evaluation equation.

\subsection{Notation}

We retain the notation \(Q_T\), \(\overline Q_T\), and \(X_T\) introduced in
the Introduction, and write \(\N\) and \(\N_0\) for the positive and
nonnegative integers, respectively.

We use standard multi-index notation for spatial derivatives. For
\(\beta=(\beta_1,\ldots,\beta_n)\in\N_0^n\), set
\[
|\beta|:=\sum_{i=1}^n\beta_i,
\qquad
D_x^\beta
:=
\partial_{x_1}^{\beta_1}\cdots\partial_{x_n}^{\beta_n}.
\]
For \(\gamma,\beta\in\N_0^n\), the relation \(\gamma\leq\beta\) is
understood componentwise.

For vector- or matrix-valued functions, \(|\cdot|\) denotes the usual
Euclidean norm on the corresponding finite-dimensional space. With this
convention, define the spatial norm
\begin{equation}\label{eq:norms-small}
    |w|_s:=\max_{|\beta|\le s}\sup_{(x,t)\in\overline Q_T}|D_x^\beta w(x,t)|,
    \qquad s\in\N_0.
\end{equation}
For $\sigma \in (0,1)$, define
\begin{equation*}
[w]_{\sigma}:=\sup_{\substack{(x,y,t)\in\R^{2n}\times[0,T]\\0<|x-y|\leq 1}}\frac{| w(x,t)- w(y,t)|}{|x-y|^\sigma}.
\end{equation*}
and, for $s\in\N_0$,
\begin{equation*}
|w|_{s,\sigma}:=|w|_s+\max_{|\beta|\le s} [D_x^\beta w]_{\sigma}.
\end{equation*}
We write $|w|_{s,0}:=|w|_s$.
These norms are purely spatial and uniform in time; they are used in
Section~\ref{sec:small-viscosity}.
For \(0<\gamma<1\), define the temporal Hölder seminorm
\[
[w]_{t;\gamma}
:=
\sup_{\substack{x\in\R^n,\ t,\tau\in[0,T]\\
0<|t-\tau|\leq1}}
\frac{|w(x,t)-w(x,\tau)|}{|t-\tau|^\gamma},
\]
and set
\[
|w|_{\alpha,\alpha/2}
:=
|w|_0+[w]_\alpha+[w]_{t;\alpha/2}.
\]
We equip
\[
X_T=C_b^{2+\alpha,1+\alpha/2}(\overline Q_T)
\]
with the standard parabolic Hölder norm
\[
\|w\|_{X_T}
:=
|w|_0
+|D_xw|_0
+|D_x^2w|_{\alpha,\alpha/2}
+|w_t|_{\alpha,\alpha/2}
+[D_xw]_{t;(1+\alpha)/2}.
\]
Finally, \(\BUC\) denotes the space of bounded uniformly continuous functions.

We also use the standard estimates: For integers \(s\geq0\) and
\(\sigma\in[0,1)\),
\begin{equation}\label{eq:tame-product}
    |vw|_{s,\sigma}\le C\bigl(|v|_{s,\sigma}|w|_0+|v|_0|w|_{s,\sigma}\bigr),
\end{equation}
and similarly for products of several factors. Moreover, if \(G\) has
bounded derivatives of order at most \(s\) (at most \(s+1\) if
\(\sigma>0\)) on \(\R^n\times B(0,r)\), and \(|w|_0\le r\), then
\begin{equation}\label{eq:tame-composition}
    |G(x,w)|_{s,\sigma}\le C\bigl(1+|w|_{s,\sigma}\bigr).
\end{equation}
Here \(C\) depends on \(n,s,\sigma\), and, in
\eqref{eq:tame-composition}, on \(r\) and the bounds on \(G\). Both
estimates follow from the Leibniz and chain rules and the interpolation
inequality \(|w|_j\le C|w|_0^{1-j/s}|w|_s^{j/s}\), \(0\le j\le s\).

\subsection{Assumptions}
\label{sec:assumptions}

We now state the structural hypotheses on the Hamiltonian $H$ in \eqref{eq:HJB}.
Throughout this section, \(T>0\) is fixed as in \ref{set:HJB}.

\begin{ass}
\label{A1}
$H$ is a $C^2$ function
with 
\[
    H, D_pH \in \BUC(\R^n\times B(0,R))
    \qquad\text{for every } R>0,
\]
and such that 
\begin{equation*}
      \lim_{|p|\to \infty} \inf_{x\in \R^n}\left(\frac{1}{2n}H(x,p)^2 + D_xH(x,p)\cdot p\right)=+\infty.
\end{equation*}
\end{ass}

\begin{rem}
Under Assumption~\ref{A1}, standard global Lipschitz estimates for
Hamilton--Jacobi equations, see \cite[Chapter 1]{tranbook} or \cite{MTL17}, imply that there
exists a constant \(C_0>0\), depending on the data but independent of \(T\),
such that, for every \(T>0\), the corresponding classical solution of
\eqref{eq:HJB} satisfies
\[
    \|u_t\|_{L^\infty(\overline Q_T)}
    +
    \|Du\|_{L^\infty(\overline Q_T)}
    \leq C_0.
\]
\end{rem}

\begin{ass}
\label{A2}
$H$ is a $C^3$ function satisfying the following bound
for any $R>0$:
\begin{equation*}
   \|H\|_{C^1(\R^n\times B(0,R))}
    +
    \|\Dp H\|_{C^1(\R^n\times B(0,R))}
    +
    \|\Dpp H\|_{C^1(\R^n\times B(0,R))}
    <\infty.
\end{equation*}
\end{ass}
\begin{rem}
\label{truncatedham}
Using a smooth cutoff in \(p\), we may replace \(H\) by a \(C^3\)
Hamiltonian that agrees with the original Hamiltonian on
\(\R^n\times B(0,2C_0+1)\) and satisfies
\begin{equation}\label{eq:conditions-H-modified}
    \|\Dp H\|_{C^1(\R^n\times \R^n)}
    +
    \|\Dpp H\|_{C^1(\R^n\times \R^n)}
    <\infty.
\end{equation}
We continue to denote the modified Hamiltonian by \(H\). Since
\(\|Du\|_{L^\infty(\overline Q_T)}\leq C_0\), this modification does not
change the equation satisfied by the exact solution \(u\).
\end{rem}

\subsection{Newton's method}\label{ss2.3}

Consider the setting of \ref{set:HJB} under
Assumptions~\ref{A1}--\ref{A2}, and let the Hamiltonian used in the iteration
be the modified Hamiltonian satisfying \eqref{eq:conditions-H-modified}.
A function \(u^0\in X_T\) is called an admissible initial guess if
\[
    u^0(x,0)=g(x)
    \qquad\text{for every }x\in\R^n.
\]

\begin{defn}[Newton iterate for the viscous Hamilton--Jacobi equation]
\label{def:newton-iterate}
Given an admissible initial guess \(u^0\), the Newton sequence
\(\{u^k\}_{k\ge0}\subset X_T\) is defined recursively as follows.  Suppose
\(u^k\in X_T\) has been constructed.  Set
\[
    b^k(x,t):=D_pH(x,Du^k(x,t)),
\]
and
\[
    c^k(x,t):=
    H(x,Du^k(x,t))
    -
    b^k(x,t)\cdot Du^k(x,t).
\]

Define \(u^{k+1}\in X_T\) to be the bounded classical solution of
\begin{equation}
\label{eq:newton-iterate-definition}
\begin{cases}
    u_t^{k+1}
    +b^k(x,t)\cdot Du^{k+1}
    +c^k(x,t)
    -\Delta u^{k+1}=0
   \qquad &\text{in }Q_T,\\[1mm]
    u^{k+1}(x,0)=g(x)
    \qquad&\text{on }\R^n.
\end{cases}
\end{equation}
Equivalently, writing
\[
    f[v]:=v_t+H(x,Dv)-\Delta v,
\]
and
\begin{equation}\label{Lk}
\mathcal L_k\varphi
    :=
    \varphi_t
    +
    D_pH(x,Du^k)\cdot D\varphi
    -
    \Delta\varphi,    
\end{equation}
the Newton correction \(h^k:=u^{k+1}-u^k\) solves
\begin{equation}
\label{eq:newton-correction-definition}
\begin{cases}
    \mathcal L_k h^k=-f[u^k]
   \qquad &\text{in }Q_T,\\[1mm]
    h^k(x,0)=0
    \qquad&\text{on }\R^n.
\end{cases}
\end{equation}
\end{defn}

\begin{rem}[Well-posedness of each Newton step]
By Assumption~\ref{A2} and \eqref{eq:conditions-H-modified} in Remark~\ref{truncatedham}, if
\(u^k\in X_T\), then the coefficients \(b^k\) and \(c^k\) in
\eqref{eq:newton-iterate-definition} belong to
\(C_b^{\alpha,\alpha/2}(\overline Q_T)\).  Hence the linear uniformly parabolic problem
\eqref{eq:newton-iterate-definition} has a unique bounded classical solution
\(u^{k+1}\in X_T\) by standard parabolic Schauder theory.
\end{rem}

The correction equation \eqref{eq:newton-correction-definition} is precisely
the Newton linearization of the PDE residual \(f\) on the affine space
\[
    \{v\in X_T:v(\cdot,0)=g\}.
\]
Indeed,
\[
    f'[u^k]\varphi
    =
    \varphi_t
    +
    D_pH(x,Du^k)\cdot D\varphi
    -
    \Delta\varphi
    =
    \mathcal L_k\varphi.
\]
Thus, recalling that \(h^k=u^{k+1}-u^k\), the equation
\[
    f'[u^k]h^k=-f[u^k]
\]
 is equivalent to
\eqref{eq:newton-iterate-definition}.

If $H$ is convex in $p$ with Legendre transform $L(x,q):=\sup_{p}\,(p\cdot q-H(x,p))$, then
\[
H(x,Du^k)-\Dp H(x,Du^k)\cdot Du^k = - L(x,D_pH(x,Du^k)),
\]
and \eqref{eq:newton-iterate-definition} is the policy-evaluation equation
associated with the corresponding supremum representation of \(H\).
The Newton iteration itself does not require convexity.

\section{Convergence of the Newton Iteration}
\label{sec:newton-convergence}

Under \eqref{eq:conditions-H-modified}, the Newton remainder admits both
global linear and local quadratic bounds in the gradient error. The linear
bound yields geometric \(L^\infty\) convergence from an arbitrary admissible
initial guess, while the quadratic bound yields local quadratic convergence
in \(X_T\).

\subsection{Global \(L^\infty\) convergence}

\begin{thm}[Global \(L^\infty\) convergence]
\label{thm:main1}
Consider the setting in \ref{set:HJB}, under Assumptions \ref{A1}--\ref{A2},
and let \(u\) be the corresponding solution.  Let \(\{u^k\}_{k\ge0}\) be the
Newton sequence of Definition \ref{def:newton-iterate}, built from the modified
Hamiltonian satisfying \eqref{eq:conditions-H-modified} and started from an
arbitrary admissible initial guess \(u^0\in X_T\).  Then, for all \(k\ge1\),
\[
    \|u^k-u\|_{L^\infty(\R^n\times[0,T])}
    \le
    C\,\|Du^0-Du\|_{L^\infty(\R^n\times[0,T])}\,2^{-k},
\]
where \(C=C(T,\|D_pH\|_{L^\infty(\R^n\times\R^n)})\).
\end{thm}
\begin{proof}[Proof of Theorem \ref{thm:main1}]
Let
\[
    v^k=u^k-u.
\]
Set
\[
    b^k(x,t)=D_pH(x,Du^k(x,t)).
\]
Subtracting \eqref{eq:HJB} from \eqref{eq:newton-iterate-definition} gives
\begin{equation}\label{eq:v-equation}
    v_t^{k+1}+b^k\cdot Dv^{k+1}+R^k=\Delta v^{k+1},
\end{equation}
where
\[
    R^k=[H(x,Du^k)-H(x,Du)]-D_pH(x,Du^k)\cdot(Du^k-Du).
\]
By the global boundedness of \(D_pH\) and the mean-value theorem, both terms
in \(R^k\) are bounded by a constant times \(|Dv^k|\).  Hence
\begin{equation}\label{eq:Rk-linear-bound}
    |R^k|\leq C|Dv^k|,
    \qquad |b^k|\leq C
    \quad\text{in }\R^n\times[0,T],
\end{equation}
where \(C\) depends only on
\(\|D_pH\|_{L^\infty(\R^n\times\R^n)}\).
Let
\[
    \phi^k=\frac12(v^k)^2.
\]
Multiplying \eqref{eq:v-equation} by $v^{k+1}$ yields
\[
    \phi_t^{k+1}+b^k\cdot D\phi^{k+1}+R^k v^{k+1}
    =\Delta\phi^{k+1}-|Dv^{k+1}|^2.
\]
Since $D\phi^{k+1}=v^{k+1}Dv^{k+1}$ and $|v^{k+1}|=\sqrt{2\phi^{k+1}}$, the
bounds \eqref{eq:Rk-linear-bound} give
$|b^k\cdot D\phi^{k+1}|\leq C|Dv^{k+1}|\sqrt{\phi^{k+1}}$ and
$|R^kv^{k+1}|\leq C|Dv^k|\sqrt{\phi^{k+1}}$.  Hence, by Young's inequality and
enlarging $C$ if necessary,
\begin{align*}
    \phi_t^{k+1}-\Delta\phi^{k+1}
    &\leq C|Dv^{k+1}|\sqrt{\phi^{k+1}}
        +C|Dv^k|\sqrt{\phi^{k+1}}-|Dv^{k+1}|^2  \\
    &\leq -\frac12|Dv^{k+1}|^2+\frac18|Dv^k|^2+C\phi^{k+1}.
\end{align*}
Hence, with
\[
    \psi^k=e^{-Ct}\phi^k,
\]
we have
\begin{equation}\label{eq:psi-ineq}
    \psi_t^{k+1}-\Delta\psi^{k+1}
    \leq -\frac12e^{-Ct}|Dv^{k+1}|^2+\frac18e^{-Ct}|Dv^k|^2.
\end{equation}
For each $k\geq0$, let $\xi^k$ be the bounded solution of
\begin{equation}\label{eq:xi-def}
\begin{cases}
    \xi_t^k-\Delta\xi^k=\dfrac12e^{-Ct}|Dv^k|^2
        \qquad&\text{in }\R^n\times(0,T),\\[2mm]
    \xi^k(x,0)=0 \qquad&\text{on }\R^n.
\end{cases}
\end{equation}
Since the source in \eqref{eq:xi-def} is nonnegative, the maximum principle for
bounded solutions of the heat equation gives $\xi^k\geq0$.
Combining \eqref{eq:psi-ineq} and \eqref{eq:xi-def} gives
\[
    \left(\psi^{k+1}+\xi^{k+1}-\frac14\xi^k\right)_t
    -\Delta\left(\psi^{k+1}+\xi^{k+1}-\frac14\xi^k\right)\leq0.
\]
Moreover, for $k\geq0$, $v^{k+1}(x,0)=0$, and $\xi^j(x,0)=0$ for all $j$, so
\[
    \left(\psi^{k+1}+\xi^{k+1}-\frac14\xi^k\right)(x,0)=0.
\]
Since \(u,u^{k+1}\in X_T\), the function \(\psi^{k+1}\) is bounded, while
\(\xi^k\) and \(\xi^{k+1}\) are bounded by construction.  Hence
\[
    w:=\psi^{k+1}+\xi^{k+1}-\frac14\xi^k
\]
is a bounded subsolution satisfying \(w(\cdot,0)=0\).  The  maximum principle therefore gives \(w\leq0\).
Since \(\psi^{k+1}\geq0\) and \(\xi^j\geq0\) for all \(j\), dropping either
term gives
\[
    0\leq \xi^{k+1}\leq \frac14\xi^k
    \qquad\text{and}\qquad
    0\leq \psi^{k+1}\leq \frac14\xi^k .
\]
Iterating the first bound yields $\xi^k\leq4^{-k}\xi^0$, and hence
$\psi^{k+1}\leq4^{-(k+1)}\xi^0$.

Finally, the function $\eta(t):=\frac12\int_0^te^{-Cs}\|Dv^0(\cdot,s)\|_{L^\infty(\R^n)}^2\,ds$
is independent of $x$ and satisfies $\eta_t-\Delta\eta\geq\frac12e^{-Ct}|Dv^0|^2$
with $\eta(0)=0$; hence $\xi^0\leq\eta$ by comparison, and
\[
    0\leq \xi^0
    \leq \frac12\int_0^t e^{-Cs}\|Dv^0(\cdot,s)\|_{L^\infty(\R^n)}^2\,ds
    \leq \frac{1}{2C}\|Dv^0\|_{L^\infty(\R^n\times[0,T])}^2 .
\]
Thus, for every \(k\geq1\),
\[
    \frac12 |v^k(x,t)|^2
    =
    \phi^k(x,t)
    =
    e^{Ct}\psi^k(x,t)
    \leq
    C_T
    \|Dv^0\|_{L^\infty(\R^n\times[0,T])}^2
    4^{-k}.
\]
Taking square roots and the supremum over \(\R^n\times[0,T]\), we obtain
\[
    \|u^k-u\|_{L^\infty(\R^n\times[0,T])}
    \leq
    C_T
    \|D(u^0-u)\|_{L^\infty(\R^n\times[0,T])}
    2^{-k},
    \qquad k\geq1.
\]
\end{proof}

\subsection{Local quadratic convergence}

The preceding proof uses only the global linear bound
\eqref{eq:Rk-linear-bound} on the Newton remainder and therefore applies to
an arbitrary admissible initial guess, at the cost of yielding only geometric
convergence in \(L^\infty\). In a neighborhood of \(u\) in \(X_T\), the
Newton remainder is quadratic in \(Du^k-Du\), and parabolic Schauder theory
converts this into a quadratic estimate in \(X_T\). As regards the
Hamiltonian, the proof below uses only the local derivative bounds in
Assumption~\ref{A2}, on a ball determined by \(\|u\|_{X_T}\).

\begin{thm}[Local quadratic convergence]
\label{thm:main2}
Assume the setting and hypotheses of Theorem~\ref{thm:main1}.
There exist constants \(\delta_0>0\) and \(C_*>0\) such that, if
\[
    \delta:=\|u^0-u\|_{X_T}\leq\delta_0,
\]
then
\[
    \|u^{k+1}-u\|_{X_T}
    \leq
    C_*\|u^k-u\|_{X_T}^2
    \qquad\text{for all }k\geq0.
\]
After decreasing \(\delta_0\) if necessary, we may assume that
\[
    \delta_0\leq1
    \qquad\text{and}\qquad
    C_*\delta_0\leq\tfrac12.
\]
Then
\[
    \|u^k-u\|_{X_T}
    \leq
    C_*^{-1}(C_*\delta)^{2^k}
    \qquad\text{for all }k\geq0.
\]
\end{thm}

\begin{proof}[Proof of Theorem \ref{thm:main2}]
Set $E_k:=\|u^k-u\|_{X_T}$, and let $b^k$ be as in Definition
\ref{def:newton-iterate}, so that \eqref{Lk} reads
$\Lop_k\varphi=\varphi_t+b^k\cdot D\varphi-\Delta\varphi$.

The integral Taylor formula in the $p$ variable gives
\[
    H(x,Du)
    =H(x,Du^k)+b^k\cdot D(u-u^k)+Q^k,
\]
where $Q^k=-R^k$ is 
%%%%XG:we already have R^k, so  Q^k is not much needed. This is a minor point.
\begin{equation}\label{eq:Qk-HJ}
    Q^k(x,t)=\int_0^1(1-s)D^2_{pp}H\bigl(x,Du^k+sD(u-u^k)\bigr)
    [D(u-u^k),D(u-u^k)]\,ds.
\end{equation}
Comparing the equation for $u$ with \eqref{eq:newton-iterate-definition} yields
\begin{equation}\label{eq:error-main2}
\begin{cases}
    \Lop_k(u^{k+1}-u)=Q^k
     \qquad&\text{in }Q_T,\\
    (u^{k+1}-u)(x,0)=0 \qquad&\text{on }\R^n.
\end{cases}
\end{equation}
Suppose $E_k\leq1$.  Then $\|u^k\|_{X_T}\leq\|u\|_{X_T}+1$, so $Du^k$ and $Du$
take values in a fixed ball $B(0,R)$ with $R=R(\|u\|_{X_T})$; by convexity, so
does $Du^k+sD(u-u^k)$ for every $s\in[0,1]$. 
By Assumption~\ref{A2} and the standard parabolic H\"older composition estimate,
\[
    \|b^k\|_{C^{\alpha,\alpha/2}(\R^n\times[0,T])}\leq C_R.
\]
Since the second-order part of $\Lop_k$ is $-\Delta$ for every $k$, the family
$\{\Lop_k\}$ is uniformly parabolic with Schauder constants independent of $k$.

Since \(u^{k+1}-u\in X_T\) is bounded and has zero initial trace, the maximum
principle applied to
\eqref{eq:error-main2} gives
\[
    \|u^{k+1}-u\|_{L^\infty(\overline Q_T)}
    \leq
    T\|Q^k\|_{L^\infty(\overline Q_T)}.
\]
The global parabolic Schauder estimate  applied to
\eqref{eq:error-main2} bounds $E_{k+1}$ by
$C\bigl(\|Q^k\|_{C^{\alpha,\alpha/2}(\overline Q_T)}
+\|u^{k+1}-u\|_{L^\infty(\overline Q_T)}\bigr)$.
Using the previous bound for the lower-order term, we obtain
\begin{equation}\label{eq:schauder-main2}
    E_{k+1}
    \leq
    C_*\|Q^k\|_{C^{\alpha,\alpha/2}(\overline Q_T)},
\end{equation}
where \(C_*\) is independent of \(k\); after the enlargement below, it
depends only on \(T,\alpha,n,\|u\|_{X_T}\) and the local \(C^1\)-bounds of
\(D_pH\) and \(D^2_{pp}H\) on \(\R^n\times B(0,R)\).

From \eqref{eq:Qk-HJ} and the algebra property of $C^{\alpha,\alpha/2}$,
\[
    \|Q^k\|_{C^{\alpha,\alpha/2}(\R^n\times[0,T])}
    \leq C_R E_k^2.
\]
Thus, after increasing $C_*$ if necessary,
\begin{equation}\label{eq:quadratic-main2}
    E_{k+1}\leq C_*E_k^2
    \qquad\text{whenever }E_k\leq1.
\end{equation}
Let \(\delta:=E_0\), and choose \(\delta_0>0\) such that
\[
    \delta_0\leq1,
    \qquad
    C_*\delta_0\leq\tfrac12.
\]
Assume that \(\delta\leq\delta_0\).  We prove by induction that
\[
    E_k\leq C_*^{-1}(C_*\delta)^{2^k}
    \qquad\text{for all }k\geq0.
\]
For \(k=0\), this is an equality.  If the estimate holds at step \(k\), then,
since \(C_*\delta\leq\tfrac12\),
\[
    E_k
    \leq C_*^{-1}(C_*\delta)^{2^k}
    \leq C_*^{-1}(C_*\delta)
    =\delta
    \leq1.
\]
Thus \eqref{eq:quadratic-main2} applies and gives
\[
    E_{k+1}
    \leq C_*E_k^2
    \leq C_*^{-1}(C_*\delta)^{2^{k+1}}.
\]
This closes the induction and proves the theorem.
\end{proof}

\section{Relation to Classical Policy Iteration}
\label{sec:policy-iteration}

We now explain how the preceding Newton iteration is related to classical
policy iteration. In particular, when \(H\) admits a sufficiently regular
control representation, the value iterates produced by the two methods
coincide. Consequently, the convergence results established above also apply
to classical policy iteration, provided the corresponding policy-evaluation
problems are well posed.

Suppose,
in addition, that the Hamiltonian used in
Definition~\ref{def:newton-iterate} admits the control representation
\begin{equation}\label{Hcontrol}
H(x,p)=\inf_{a\in A}\left[c(x,a)+p\cdot b(x,a)\right].
\end{equation}
Here we assume that \(A\) is a compact metric space and that \(c\) and \(b\)
are continuous and bounded, so that the infimum in \eqref{Hcontrol} is
attained.

Note that \eqref{Hcontrol} forces \(H(x,\cdot)\) to be concave, with
\(|D_pH|\le\|b\|_{L^\infty}\). Since the proof of Theorem~\ref{thm:main1}
uses only a global bound on \(D_pH\), the cutoff of
Remark~\ref{truncatedham} is not needed here, and \eqref{Hcontrol} may be
the representation of the original Hamiltonian. For a
supremum representation of a convex Hamiltonian, the analogous discussion
applies with \(\inf\) and \(\argmin\) replaced by \(\sup\) and \(\argmax\),
respectively.
We formulate the iteration under the proviso that the required minimizing
selectors can be chosen so that the corresponding policy-evaluation problems
admit unique solutions in \(X_T\).
\begin{defn}[Policy Iteration]
Given an initial policy \(a^0(x,t)\), for \(k\in\N_0\), let \(u^k\) be
the bounded classical solution of
\begin{equation*}
\begin{cases}
\partial_t u^k + c(x,a^k(x,t)) + D u^k \cdot b(x,a^k(x,t)) -\Delta u^k=0  \qquad &\text{in } Q_T,\\
u^k(x,0) = g(x) \qquad &\text{on }  \mathbb{R}^n. 
 \end{cases}
 \end{equation*}
Then choose a minimizing selector
\begin{equation}
\label{eq:uiter}
a^{k+1}(x,t)\in  \argmin_{a\in A} \left[c(x,a) + D u^k(x,t) \cdot b(x,a)\right].
\end{equation}
\end{defn}
Although the boundedness of \(b\) makes the Hamiltonian in
\eqref{Hcontrol} Lipschitz in \(p\), 
the policy update
\eqref{eq:uiter} requires \(Du^k\) to be defined pointwise and a minimizing
selector whose induced coefficients yield a well-posed policy-evaluation
problem. In general control problems, these requirements may fail, and the
undiscretized policy iteration may be ill-posed; see \cite{TWZ}. By contrast, Newton's method does not involve the
controls directly, but it requires suitable smoothness of \(H\).

The following lemma says that the two iteration schemes are essentially the same.
\begin{lem}
Fix \(k\geq0\), and suppose that \(H(x,\cdot)\) is differentiable at
\(p=Du^k(x,t)\) for every \((x,t)\in\overline Q_T\). Then, for every
\((x,t)\in\overline Q_T\),
\[
c^k(x,t)+b^k(x,t)\cdot Du^{k+1}
    =c(x,a^{k+1}(x,t)) +  b(x,a^{k+1}(x,t))\cdot D u^{k+1}
\]
where $b^k$ and $c^k$ are given in Definition \ref{def:newton-iterate}.
Consequently, whenever the two schemes have the same value iterate \(u^k\),
they determine the same linear problem for \(u^{k+1}\). In particular, if
the value function \(u^0\) produced from the initial policy \(a^0\) is used
as the initial Newton guess, then, provided all steps are well posed, the two
value-function sequences coincide.
\end{lem}

\begin{proof}
Fix \((x,t)\), and set \(p=Du^k(x,t)\). By \eqref{Hcontrol} and
\eqref{eq:uiter}, the affine function
\[
    q\longmapsto
    c(x,a^{k+1}(x,t))+q\cdot b(x,a^{k+1}(x,t))
\]
touches the concave function \(H(x,\cdot)\) from above at \(p\).
Therefore \(b(x,a^{k+1}(x,t))\) is a supergradient of \(H(x,\cdot)\) at
\(p\). Since \(H(x,\cdot)\) is differentiable at \(p\),
\[
    D_pH(x,Du^k(x,t))
    =
    b(x,a^{k+1}(x,t)).
\]
The definition of $a^{k+1}$ yields
\[
H(x,Du^k(x,t))=c(x,a^{k+1}(x,t)) + D u^k \cdot b(x,a^{k+1}(x,t)).
\]
Thus, also using the definition of $b^k$ and $c^k$, we have
\begin{align*}
&\quad\,\, c^k(x,t)+b^k(x,t)\cdot Du^{k+1} = H(x,Du^k)
    +
    D_pH(x,Du^k)\cdot ( Du^{k+1}-Du^k)\\
    &=c(x,a^{k+1}(x,t)) + D u^k \cdot b(x,a^{k+1}(x,t))+D_pH(x,Du^k)\cdot ( Du^{k+1}-Du^k)\\
    &=c(x,a^{k+1}(x,t)) +  b(x,a^{k+1}(x,t))\cdot D u^{k+1}.    
\end{align*}
Thus, the policy-evaluation equation for \(u^{k+1}\) is exactly the Newton
equation linearized at \(u^k\). Uniqueness of the corresponding linear
problems and induction prove the final assertion.
\end{proof}

\begin{cor}[Convergence of classical policy iteration]
\label{cor:policy-convergence}
Assume the setting and hypotheses of Theorem~\ref{thm:main1},
the control representation
\eqref{Hcontrol}, and that the policy-evaluation problem for \(a^0\) has a
unique bounded classical solution \(u^0\in X_T\).
Suppose that the value function \(u^0\) produced
from the initial policy \(a^0\) is used as the initial Newton guess. Then the
value-function iterates of classical policy iteration coincide with the
Newton iterates of Definition~\ref{def:newton-iterate}. 
Consequently, they satisfy the global convergence estimate of
Theorem~\ref{thm:main1} and, whenever the initial-value closeness condition
of Theorem~\ref{thm:main2} holds, its local quadratic convergence estimate.
\end{cor}
\begin{proof}
The proof of the preceding lemma shows that \(b(x,a^{k+1})=b^k\) and
\(c(x,a^{k+1})=c^k\) for every minimizing selector. Hence each subsequent
policy-evaluation problem coincides with a well-posed Newton step, and the
claim follows from Theorems~\ref{thm:main1}--\ref{thm:main2}.
\end{proof}

\section{Fully Nonlinear Parabolic Equations}
\label{sec:fully-nonlinear}

The local quadratic convergence mechanism of
Theorem~\ref{thm:main2} extends to uniformly parabolic fully nonlinear
equations. Fix \(T>0\) and consider
\begin{equation}\label{fully-nonlinear}
\begin{cases}
u_t+F(x,Du,D^2u)=0
\qquad&\text{in }Q_T,\\
u(x,0)=g(x)
\qquad&\text{on }\R^n.
\end{cases}
\end{equation}
In the fully nonlinear setting, the regularity assumptions on the
Hamiltonian and the fixed linear diffusion \(-\Delta u\) are replaced by
the corresponding regularity and uniform ellipticity conditions on \(F\).

Let \(\cS(n)\) be the space of real symmetric \(n\times n\) matrices, and write \(F=F(x,p,X)\) for
\[
(x,p,X)\in\R^n\times\R^n\times\cS(n).
\]

For \(A,B\in\cS(n)\), we define the Frobenius inner product
\[
A:B:=\operatorname{tr}(AB).
\]
We use \(D_pF\) and \(D_XF\) to denote 
differentiation with respect to
the second and third variables, respectively. We identify \(D_XF\) as an element of \(\cS(n)\) through the Frobenius inner product.
For \(R>0\), set
\[
\mathcal D_R
:=
\R^n\times B(0,R)\times B_{\cS(n)}(0,R).
\]

\begin{ass}[Regularity and uniform ellipticity]
\label{ass:fully-nonlinear}
% We write \(F=F(x,p,X)\), where
% \[
% (x,p,X)\in\R^n\times\R^n\times\cS(n).
% \]
% The symbols \(D_pF\) and \(D_XF\) denote differentiation with respect to
% the second and third variables, respectively, with \(D_XF\) identified
% with an element of \(\cS(n)\) through the Frobenius inner product.
% For \(R>0\), set
% \[
% \mathcal D_R
% :=
% \R^n\times B(0,R)\times B_{\cS(n)}(0,R).
% \]
%%%%XG:the sentences above are definitions (not parts of the assumption) which makes Assumption 3 too long, so I moved them to the above.
Let \(F\in C^3(\R^n\times\R^n\times\cS(n))\). For every \(R>0\),
\begin{align}
&|F|_{C^1(\mathcal D_R)}
+|D_pF|_{C^1(\mathcal D_R)}
+|D_XF|_{C^1(\mathcal D_R)}
\notag\\
&\qquad
+|D_{pp}^2F|_{C^1(\mathcal D_R)}
+|D_{pX}^2F|_{C^1(\mathcal D_R)}
+|D_{XX}^2F|_{C^1(\mathcal D_R)}
<\infty.
\label{eq:F-regularity}
\end{align}
Moreover, there exist \(0<\lambda\leq\Lambda\) such that, for every
\((x,p,X)\in\R^n\times\R^n\times\cS(n)\) and every
\(Y\in\cS(n)\) with \(Y\geq0\),
\begin{equation}
\label{eq:F-ellipticity}
-\Lambda\operatorname{tr}(Y)
\leq F(x,p,X+Y)-F(x,p,X)
\leq-\lambda\operatorname{tr}(Y).
\end{equation}
\end{ass}

Let \(u^0\in X_T\) satisfy \(u^0(\cdot,0)=g\). We define the Newton
sequence recursively. For \(v\in X_T\), write
\[
F[v]:=F(x,Dv,D^2v),
\]
and use the same evaluation convention for the derivatives of \(F\).
For \(k\geq0\), set
\[
\ell^k(x,t)
:=
F[u^k]
-D_pF[u^k]\cdot Du^k
-D_XF[u^k]:D^2u^k.
\]
Given \(u^k\in X_T\), let \(u^{k+1}\in X_T\) be the unique bounded
classical solution of
\begin{equation}\label{eq:iter-F}
\begin{cases}
u_t^{k+1}
+D_pF[u^k]\cdot Du^{k+1}
+D_XF[u^k]:D^2u^{k+1}
+\ell^k(x,t)=0
\qquad&\text{in }Q_T,\\
u^{k+1}(x,0)=g(x)
\qquad&\text{on }\R^n.
\end{cases}
\end{equation}
Indeed, Assumption~\ref{ass:fully-nonlinear} implies that
\(D_pF[u^k]\), \(D_XF[u^k]\), and \(\ell^k\) belong to
\(C_b^{\alpha,\alpha/2}(\overline Q_T)\), 
while
\(-D_XF[u^k]\) is uniformly elliptic: applying
\eqref{eq:F-ellipticity} with \(\tau Y\), dividing by \(\tau>0\), and
letting \(\tau\downarrow0\) gives \(\lambda I\leq-D_XF\leq\Lambda I\) on
\(\R^n\times\R^n\times\cS(n)\).

Thus each Newton step is well posed in \(X_T\) by standard parabolic
Schauder theory.

\begin{thm}[Fully nonlinear local quadratic convergence]
\label{thm:main3}
Suppose that Assumption~\ref{ass:fully-nonlinear} holds. Let
\((u^k)_{k\geq0}\subset X_T\) be the Newton sequence defined by
\eqref{eq:iter-F}, and let \(u\in X_T\) be a classical solution of
\eqref{fully-nonlinear}. Then there exist constants
\(C_*>0\) and \(\delta_0>0\) such that, if
\[
\delta:=\|u^0-u\|_{X_T}\leq\delta_0,
\]
then
\[
\|u^{k+1}-u\|_{X_T}
\leq
C_*\|u^k-u\|_{X_T}^2
\qquad\text{for every }k\geq0.
\]
Moreover, \(\delta_0\) may be chosen so that
\[
\delta_0\leq1,
\qquad
C_*\delta_0\leq\tfrac12,
\]
and then
\[
\|u^k-u\|_{X_T}
\leq
C_*^{-1}(C_*\delta)^{2^k}
\qquad\text{for every }k\geq0.
\]
\end{thm}
\begin{proof}
Set
\[
E_k:=\|u^k-u\|_{X_T},
\qquad
P^k:=Du^k,\quad X^k:=D^2u^k,
\qquad
P:=Du,\quad X:=D^2u.
\]
Define the linearized operator
\[
    \Lop_k\varphi
    =\varphi_t+D_pF[u^k]\cdot D\varphi+D_XF[u^k]:D^2\varphi.
\]
Since \(\lambda I\leq-D_XF\leq\Lambda I\), as noted before the theorem,
\(\Lop_k\) is uniformly parabolic.

By the integral Taylor formula in the variables $(p,X)$,
\begin{align*}
    F(x,Du,D^2u)
    &=F[u^k]+D_pF[u^k]\cdot D(u-u^k)
        +D_XF[u^k]:D^2(u-u^k)+r^k,
\end{align*}
where
\begin{align}\label{eq:rk-F}
    r^k
    =\int_0^1(1-s)&\Big(
    D^2_{pp}F(x,P^k+s(P-P^k),X^k+s(X-X^k))[P-P^k,P-P^k]
    \\\notag &+2D^2_{pX}F(x,P^k+s(P-P^k),X^k+s(X-X^k))[P-P^k,X-X^k]
    \\\notag &+D^2_{XX}F(x,P^k+s(P-P^k),X^k+s(X-X^k))[X-X^k,X-X^k]
    \Big)\,ds.
\end{align}
Comparing this identity with the definition of $u^{k+1}$ gives
\begin{equation}\label{eq:error-main3}
\begin{cases}
    \Lop_k(u^{k+1}-u)=r^k
        \qquad&\text{in }Q_T,\\[1mm]
    (u^{k+1}-u)(x,0)=0 \qquad&\text{on }\R^n.
\end{cases}
\end{equation}
Assume \(E_k\leq1\). Then
\[
\|u^k\|_{X_T}\leq \|u\|_{X_T}+1.
\]
Consequently, for some \(R=R(\|u\|_{X_T})\), all the arguments of \(F\)
along the segment in \eqref{eq:rk-F} belong to \(\mathcal D_R\).
By \eqref{eq:F-regularity}, the coefficients \(D_pF[u^k]\) and
\(D_XF[u^k]\) have uniformly bounded
\(C_b^{\alpha,\alpha/2}(\overline Q_T)\) norms, with bounds independent
of \(k\).
Since \(u^{k+1}-u\) has zero initial trace, the maximum principle applied
to \eqref{eq:error-main3} gives
\[
|u^{k+1}-u|_0
\leq
T|r^k|_0.
\]
Combining this estimate with the global parabolic Schauder estimate on
\(\R^n\times[0,T]\) yields
\[
E_{k+1}
\leq
C|r^k|_{\alpha,\alpha/2},
\]
where \(C\) is independent of \(k\).
Using \eqref{eq:rk-F}, the \(C^1\) bounds on the second derivatives of
\(F\) in \eqref{eq:F-regularity}, and the standard parabolic Hölder
composition and product estimates, we obtain
\[
|r^k|_{\alpha,\alpha/2}
\leq
C_RE_k^2,
\]
where \(C_R\) is independent of \(k\). Combining the preceding two
estimates and enlarging the resulting constant if necessary, we obtain
\begin{equation}\label{eq:quadratic-main3}
    E_{k+1}\leq C_*E_k^2
    \qquad\text{whenever }E_k\leq1.
\end{equation}
Choose \(\delta_0>0\) so that
\[
\delta_0\leq1,
\qquad
C_*\delta_0\leq\tfrac12.
\]
Let \(\delta:=E_0\leq\delta_0\). We claim that
\[
E_k
\leq
C_*^{-1}(C_*\delta)^{2^k}
\qquad\text{for every }k\geq0.
\]
The claim is an equality for \(k=0\). If it holds at step \(k\), then,
since \(C_*\delta\leq\tfrac12\),
\[
E_k
\leq
C_*^{-1}(C_*\delta)^{2^k}
\leq
C_*^{-1}(C_*\delta)
=
\delta
\leq1,
\]
and hence \eqref{eq:quadratic-main3} gives
\[
E_{k+1}
\leq
C_*E_k^2
\leq
C_*^{-1}(C_*\delta)^{2^{k+1}}.
\]
This closes the induction.
\end{proof}

\section{A Small-Viscosity Nash--Moser Newton Iteration}
\label{sec:small-viscosity}

The preceding local convergence arguments use parabolic Schauder
estimates that are not uniform as the viscosity tends to zero. We
therefore modify the Newton step by smoothing both the current iterate
and its residual in space, while preserving the prescribed initial
condition. The analysis combines transport estimates uniform in the
viscosity with spatial H\"older estimates whose viscosity dependence is
explicit. Theorem~\ref{thm:criterion-new} gives a conditional convergence
criterion for this Nash--Moser scheme: under a smallness condition on the
initial residual compatible with the viscosity and the smoothing scale,
classical iterates have decaying residuals and converge in the spatial
norm $|\cdot|_{m+1}$.

\subsection{The iteration and convergence criterion}

We work directly with $H$ under the hypotheses stated below, without
using the cutoff of Remark~\ref{truncatedham}. For $0<\ep\le1$,
consider the small-viscosity Hamilton--Jacobi equation
\begin{equation}\label{eq:viscousHJ-eps}
\begin{cases}
     u_t+H(x,Du)-\ep\Delta u=0
     \qquad&\text{in } Q_T,\\   
    u(x,0)=g(x) \qquad&\text{on } \R^n.
\end{cases}
\end{equation}
Set
\begin{equation}\label{eq:Feps}
    f_\ep[u]=u_t+H(x,Du)-\ep\Delta u,
\end{equation}
and
\begin{equation}\label{eq:linearized-eps}
    \Lop_v^\ep h=h_t+D_pH(x,Dv)\cdot Dh-\ep\Delta h.
\end{equation}
Let \(N_*\) be a positive integer. Choose a kernel
 $\rho\in C_c^\infty(\R^n)$ satisfying
\begin{equation}\label{eq:kernel-moments}
    \int_{\R^n}\rho(z)\,dz=1,
    \qquad
    \int_{\R^n}z^\beta\rho(z)\,dz=0
    \quad \text{for }1\le |\beta|\le N_*.
\end{equation}
For background on smoothing operators in Nash--Moser schemes, see
\cite{hamilton1982inverse}.
When $N_*\ge2$, the moment conditions force $\rho$ to change sign.
% This cancellation is
The vanishing moments in \eqref{eq:kernel-moments}
will be used in the high-order
approximation estimates below.
For $\theta\ge1$, set
\begin{equation}\label{eq:convolution-smoothing}
    \rho_\theta(x)=\theta^n\rho(\theta x),
    \qquad
    S_\theta w=\rho_\theta *_x w,
\end{equation}
where $*_x$ means that the convolution is only in the spatial variable.

Given an initial guess \(u^0:\overline Q_T\to\R\) satisfying
the regularity assumptions below and \(u^0(\cdot,0)=g\),
we define \((u^k)_{k\geq0}\) recursively.
Set
\begin{equation}\label{eq:residual-new}
    R_k=f_{\ep}[u^k].
\end{equation}
Let \(\bar g(x,t):=g(x)\) for \((x,t)\in\overline Q_T\). For a smoothing
parameter \(\theta_k\geq1\), define the smoothed iterate
\begin{equation}\label{eq:smoothed-iterate-new}
\widetilde u^k:=\bar g+S_{\theta_k}(u^k-\bar g).
\end{equation}
Since \(u^k(\cdot,0)=g\) by induction, \(u^k-\bar g\) has zero initial
trace and hence \(\widetilde u^k(\cdot,0)=g\). The smoothed Newton step is
\begin{equation}\label{eq:NMstep-new}
\begin{cases}
    \Lop_{\widetilde u^k}^{\ep} h^k=-S_{\theta_k}R_k
        \qquad&\text{in } Q_T,\\
    h^k(x,0)=0 \qquad&\text{on }\R^n,
\end{cases}
\end{equation}
and
\begin{equation}\label{eq:update-new}
    u^{k+1}=u^k+h^k.
\end{equation}

\begin{thm}[Small-viscosity Nash--Moser convergence criterion]\label{thm:criterion-new}
Let \(T>0\), let \(m\geq1\) be an integer, and let \(\mu\in[0,1)\).
Assume that the kernel
in \eqref{eq:kernel-moments} is chosen with
\[
N_*\geq m+3.
\]
Assume that
\[
H\in C_b^{m+6}(\R^n\times B(0,r))
\qquad\text{for every }r>0,
\]
and that
\[
u^0\in C_b^{m+6}(\overline Q_T),
\qquad u^0(\cdot,0)=g.
\]
Then there exist \(a,\chi\in(1,2)\), depending only on \(\mu\), and a
threshold \(\Theta_*\geq2\) such that the following holds. Fix
\(\theta_0\geq\Theta_*\) and set
\begin{equation*}
\theta_k:=\theta_0^{\chi^k},
\qquad k\geq0.
\end{equation*}
Assume that
\begin{equation}\label{eq:bootstrap-new'}
\theta_0^{-\mu}\leq\ep\leq1,
\end{equation}
and
\begin{equation}\label{eq:initial-small-new}
|R_0|_m\leq\theta_0^{-a}.
\end{equation}
Then the iteration \eqref{eq:NMstep-new}--\eqref{eq:update-new}
is well defined for every $k\ge0$, and its iterates satisfy
\begin{equation}\label{eq:residual-decay-new}
    |R_k|_m\le \theta_0^{-a\chi^k}
    \qquad\text{for all }k\ge0.
\end{equation}
Moreover, \((u^k)\) converges with respect to \(|\cdot|_{m+1}\) to a
bounded classical solution \(u\) of \eqref{eq:viscousHJ-eps}, and
\begin{equation*}
|u^k-u|_{m+1}
\leq
C\bigl(\theta_0^{1-a}\bigr)^{\chi^k},
\qquad k\geq0.
\end{equation*}
Here \(\Theta_*\) and \(C\) depend on
\(m,\mu,T,g,H,u^0\), and the fixed kernel \(\rho\), but
are independent of \(0<\ep\leq1\) and \(\theta_0\).
\end{thm}

The proof combines spatial smoothing estimates, linear estimates that
track the dependence on \(\ep\), and a one-step estimate for the nonlinear
residual. We establish these ingredients in the next three subsections
before proving Theorem~\ref{thm:criterion-new}.

\subsection{Smoothing estimates}

Recall $S_\theta$ from \eqref{eq:convolution-smoothing}.
The first lemma concerns the $C^{s,\sigma}$ regularity of $S_\theta$.

\begin{lem}\label{lem:smoothing}
For integers $0\le s\le r$ and $\sigma\in[0,1)$, there are
constants $C_s$ and $C_{r,s,\sigma}$, depending only on the
indicated indices, $n$, and the fixed kernel $\rho$, such that,
for all $\theta\ge1$ and every $w$ with bounded continuous
spatial derivatives through order $s$, whenever the norm on
the right-hand side of the corresponding estimate is finite,
\begin{align}
    |S_\theta w|_{s,\sigma}&\le C_{s}|w|_{s,\sigma},\label{eq:smoothing1-new}\\
    |S_\theta w|_{r,\sigma}&\leq C_{r,s,\sigma}\theta^{r+\sigma-s} |w|_s.\label{eq:smoothing2-new}
\end{align}
\end{lem}
\begin{proof}
Since $D^\beta S_\theta w=\rho_\theta*D^\beta w$ and
$\|\rho_\theta\|_{L^1}=\|\rho\|_{L^1}$, 
convolution with $\rho_\theta$ increases the supremum norm and,
for $0<\sigma<1$, the H\"older seminorm by at most the factor
$\|\rho\|_{L^1}$. This gives \eqref{eq:smoothing1-new} with
$C_s=\|\rho\|_{L^1}$.
For \eqref{eq:smoothing2-new}, we place at most $s$ derivatives on $w$ and
the remaining ones on the kernel, where each derivative costs a factor
$\theta$. Fix a multi-index $\beta$ with
$|\beta|\leq r$, and choose $\gamma\leq\beta$ such that
\[
|\gamma|=\min\{|\beta|,s\}.
\]
Set $q:=|\beta|-|\gamma|=(|\beta|-s)_+$. Then
\[
D^\beta S_\theta w=(D^{\beta-\gamma}\rho_\theta)*D^\gamma w,
\]
and, since
$\|D^{\beta-\gamma}\rho_\theta\|_{L^1}
=\theta^q\|D^{\beta-\gamma}\rho\|_{L^1}$,
\[
|D^\beta S_\theta w|_0
\leq \|D^{\beta-\gamma}\rho_\theta\|_{L^1}|w|_s
\leq C\theta^q|w|_s
\leq C\theta^{r-s}|w|_s.
\]
Taking the maximum over $|\beta|\leq r$ proves
\eqref{eq:smoothing2-new} when $\sigma=0$.

Suppose now that $0<\sigma<1$. It remains to bound
$[D^\beta S_\theta w]_\sigma$; here a fraction $\sigma$ of a derivative
costs a factor $\theta^\sigma$. Set
$\widetilde\rho:=D^{\beta-\gamma}\rho$.
Fix $t\in[0,T]$ and $x,y\in\R^n$ with $0<|x-y|\leq1$, and let
$\zeta:=\theta(x-y)$. For $v\in\R^n$, write
$\tau_\zeta\widetilde\rho(v):=\widetilde\rho(v+\zeta)$.
Then
\[
D^\beta S_\theta w(x,t)-D^\beta S_\theta w(y,t)
=\int_{\R^n}\bigl(D^{\beta-\gamma}\rho_\theta(x-z)
-D^{\beta-\gamma}\rho_\theta(y-z)\bigr)D^\gamma w(z,t)\,dz.
\]
Since
$D^{\beta-\gamma}\rho_\theta(v)
=\theta^{n+q}\widetilde\rho(\theta v)$, the change of variables
$v=\theta(y-z)$, for which $x-z=(v+\zeta)/\theta$ and
$dz=\theta^{-n}dv$, gives
\[
\begin{aligned}
&\int_{\R^n}
\bigl|D^{\beta-\gamma}\rho_\theta(x-z)
      -D^{\beta-\gamma}\rho_\theta(y-z)\bigr|\,dz \\
&\qquad
=\theta^q\|\tau_\zeta\widetilde\rho-\widetilde\rho\|_{L^1(\R^n)}.
\end{aligned}
\]
On the one hand,
$\|\tau_\zeta\widetilde\rho-\widetilde\rho\|_{L^1}
\leq2\|\widetilde\rho\|_{L^1}$;
on the other hand, by the mean value theorem,
$\|\tau_\zeta\widetilde\rho-\widetilde\rho\|_{L^1}
\leq|\zeta|\,\|D\widetilde\rho\|_{L^1}$.
Combining these bounds with
$\min\{1,|\zeta|\}\le|\zeta|^\sigma$ gives
\[
\|\tau_\zeta\widetilde\rho-\widetilde\rho\|_{L^1}
\leq C|\zeta|^\sigma,
\qquad
C:=\max\bigl\{
2\|\widetilde\rho\|_{L^1},
\|D\widetilde\rho\|_{L^1}
\bigr\}.
\]
Here $C$ depends on $\beta$ and $\gamma$, of which there are finitely
many. Consequently, by the identity above and $\zeta=\theta(x-y)$,
\[
\begin{aligned}
\bigl|D^\beta S_\theta w(x,t)-D^\beta S_\theta w(y,t)\bigr|
&\leq
\theta^q\|\tau_\zeta\widetilde\rho-\widetilde\rho\|_{L^1}
\,|D^\gamma w|_0 \\
&\leq
C\theta^{q+\sigma}|x-y|^\sigma|w|_s,
\end{aligned}
\]
so that
\[
[D^\beta S_\theta w]_\sigma
\leq C\theta^{q+\sigma}|w|_s
\leq C\theta^{r+\sigma-s}|w|_s.
\]
Taking the maximum over $|\beta|\leq r$ and combining with the
supremum-norm estimate above, using $\theta\geq1$, proves
\eqref{eq:smoothing2-new}.
\end{proof}

The next lemma bounds $(I-S_\theta)w$ in terms of the higher norm of $w$. It relies on the cancellation property \eqref{eq:kernel-moments}.

\begin{lem}\label{lem:smoothing'}
Under the assumptions of Lemma \ref{lem:smoothing}, further assume that $r+\sigma\leq N_*+1$. 
Then there is a constant $C_{r,s,\sigma}>0$, depending only on
$n,r,s,\sigma$, and the fixed kernel $\rho$, such that, for all
$\theta\ge1$ and every $w$ with bounded continuous spatial
derivatives through order $r$ and $|w|_{r,\sigma}<\infty$,
\begin{align}
    |(I-S_\theta)w|_{s}&\le C_{r,s,\sigma}\theta^{s-r-\sigma}|w|_{r,\sigma}.\label{eq:smoothing3-new}
\end{align}
\end{lem}
\begin{proof}
Fix $t\in[0,T]$ and $|\beta|\le s$, and set
\[
v(x):=D_x^\beta w(x,t),
\qquad k:=r-|\beta|.
\]
By the definition of convolution,
\[
D_x^\beta S_\theta w(x,t)
=\int_{\R^n}\rho_\theta(z)v(x-z)\,dz.
\]

If $k=0$, then
\[
D_x^\beta(I-S_\theta)w(x,t)
=\int_{\R^n}\rho_\theta(z)\bigl(v(x)-v(x-z)\bigr)\,dz.
\]
Using boundedness when $\sigma=0$, and, when $0<\sigma<1$, H\"older
continuity for $|z|\le1$ and boundedness for $|z|>1$, we obtain
\[
|D_x^\beta(I-S_\theta)w(x,t)|
\le C\theta^{-\sigma}|w|_{r,\sigma}.
\]

Suppose now that $k\ge1$. Taylor's formula gives
\begin{equation}\label{5.4}
\begin{aligned}
v(x-z)
={}&\sum_{|\alpha|<k}
\frac{(-z)^\alpha}{\alpha!}D^\alpha v(x)\\
&+k\sum_{|\alpha|=k}\frac{(-z)^\alpha}{\alpha!}
\int_0^1(1-\tau)^{k-1}D^\alpha v(x-\tau z)\,d\tau.
\end{aligned}
\end{equation}
For $\sigma=0$, the constant Taylor term cancels by normalization,
and the moments of degrees $1,\ldots,k-1$ vanish by \eqref{eq:kernel-moments}, since
$k-1\le r-1\le N_*$. Hence
\[
|D_x^\beta(I-S_\theta)w(x,t)|
\le C|w|_r\int_{\R^n}|z|^k|\rho_\theta(z)|\,dz
\le C\theta^{-k}|w|_r.
\]

For $0<\sigma<1$, the assumption $r+\sigma\le N_*+1$ implies
$r\le N_*$. Thus the moments of degree $k$ also vanish, and we
may subtract $D^\alpha v(x)$ inside the remainder integral.
Moreover,
\[
|D^\alpha v(x-\tau z)-D^\alpha v(x)|
\le C|w|_{r,\sigma}\tau^\sigma|z|^\sigma
\qquad (|\alpha|=k).
\]
For $|\tau z|\le1$ this follows from the H\"older seminorm, and
for $|\tau z|>1$ it follows from the supremum norm. Consequently,
\[
\begin{aligned}
|D_x^\beta(I-S_\theta)w(x,t)|
&\le C|w|_{r,\sigma}
\int_{\R^n}|z|^{k+\sigma}|\rho_\theta(z)|\,dz\\
&\le C\theta^{-k-\sigma}|w|_{r,\sigma}.
\end{aligned}
\]
Taking the supremum over $(x,t)$ and the maximum over
$|\beta|\le s$, and using $\theta\ge1$, proves
\eqref{eq:smoothing3-new}.
\end{proof}

\subsection{Linear estimates}

We first derive integer-order estimates uniform in $\ep$. A heat-kernel
estimate then provides spatial H\"older control with explicit dependence
on $\ep$.

\begin{lem}\label{lem:linear-new}
Let $s\ge0$ be an integer. Suppose that
$b:\overline Q_T\to\R^n$ and $f:\overline Q_T\to\R$
are continuous and have bounded continuous spatial derivatives
through order $\max\{s,1\}$. Let $h$ be a bounded classical
solution of
\begin{equation}\label{eq:linear-advection-diffusion-new}
\begin{cases}
    h_t+b(x,t)\cdot Dh-\ep\Delta h=f,\\
    h(x,0)=0.
\end{cases}
\end{equation}
Assume $|b|_1\le M$.  Then
\begin{equation}\label{eq:linear-estimate-new}
    |h|_s\le C_s(M,T)\bigl(|f|_s+|b|_s|f|_1\bigr),
\end{equation}
where $C_s(M,T)$ is independent of $0<\ep\le1$.
\end{lem}
\begin{proof}
For \(s=0\), the maximum principle gives \(|h|_0\le T|f|_0\),
which implies \eqref{eq:linear-estimate-new}.
Let \(s\ge1\), and assume first that \(b,f\) are smooth in \(x\)
with bounded spatial derivatives of all orders. For fixed \(\ep>0\),
spatial parabolic regularity gives bounded spatial derivatives of
\(h\) of all orders, so we may differentiate the equation and apply
the maximum principle below.

For \(|\beta|\le s\), differentiating
\eqref{eq:linear-advection-diffusion-new} yields
\begin{equation}\label{eq:commutator-new}
    (D^\beta h)_t+b\cdot D(D^\beta h)-\ep\Delta(D^\beta h)
    =D^\beta f-[D^\beta,b\cdot D]h.
\end{equation}
For \(t\in[0,T]\) and \(j\le s\), write
\[
    |h(\cdot,t)|_j
    :=\max_{|\gamma|\le j}\sup_{x\in\R^n}|D^\gamma h(x,t)|.
\]
Leibniz's rule, spatial interpolation, and \(|b|_1\le M\) give
\begin{equation}\label{eq:commutator-estimate-new}
    |[D^\beta,b\cdot D]h(\cdot,t)|_0
    \le C_s\bigl(M|h(\cdot,t)|_s
                 +|b|_s|h(\cdot,t)|_1\bigr),
\end{equation}
where \(C_s\) depends only on \(n\) and \(s\).
Since \(D^\beta h(\cdot,0)=0\), the maximum principle yields
\[
    |h(\cdot,t)|_s
    \le t|f|_s
    +C_s\int_0^t
    \bigl(M|h(\cdot,\tau)|_s
          +|b|_s|h(\cdot,\tau)|_1\bigr)\,d\tau.
\]
For \(s=1\), Gronwall's inequality and \(|b|_1\le M\) give
\[
    |h|_1\le C_1(M,T)|f|_1.
\]
For \(s\ge2\), substituting this bound into the preceding inequality
and applying Gronwall's inequality again proves
\eqref{eq:linear-estimate-new}, with a constant independent of \(\ep\).

Finally, let \(b,f\) be as in the statement, and let
\(b^\eta,f^\eta\) be their spatial convolutions with a nonnegative
unit-mass mollifier. Then
\[
    |b^\eta|_j\le|b|_j,
    \qquad |f^\eta|_j\le|f|_j,
    \qquad 0\le j\le s.
\]
In particular, \(|b^\eta|_1\le M\). The spatial Lipschitz bounds
also give \(b^\eta\to b\) and \(f^\eta\to f\) uniformly.
Let \(h^\eta\) be the bounded solution of
\eqref{eq:linear-advection-diffusion-new} with data \(b^\eta,f^\eta\).
By the smooth case, \(h^\eta\) satisfies
\eqref{eq:linear-estimate-new} and
\(|h^\eta|_1\le C_1(M,T)|f|_1\), uniformly in \(\eta\).
Comparison gives
\[
    |h^\eta-h|_0
    \le T\bigl(|f^\eta-f|_0
              +|b^\eta-b|_0\,|Dh^\eta|_0\bigr)
    \longrightarrow0.
\]
For \(|\beta|\le s-1\), \(D^\beta h^\eta\) has zero initial data
and solves the heat equation with source
\(D^\beta(f^\eta-b^\eta\cdot Dh^\eta)\), bounded uniformly in \(\eta\).
For fixed \(\ep>0\), Lemma~\ref{L.7.6}, whose proof is independent
of the present lemma, therefore gives uniform spatial H\"older
bounds through order \(s\).
At each fixed \(t\), the Arzel\`a--Ascoli theorem yields locally uniform
convergence of these derivatives along a subsequence.
Since \(h^\eta\to h\) uniformly, their limits are the spatial
derivatives of \(h(\cdot,t)\).
Passing the integer derivative bounds to the limit and taking
the supremum over \(t\in[0,T]\) proves
\eqref{eq:linear-estimate-new} with the same
\(\ep\)-independent constant.
\end{proof}

The following lemma concerns the spatial H\"{o}lder norm of $Dh(\cdot,t)$.

\begin{lem}\label{L.7.6}
Let $\Gamma:\overline Q_T\to\R$ be bounded and measurable,
and let $h$ be the bounded mild solution of
$$
\begin{cases} 
h_t - \varepsilon \Delta h = \Gamma(x,t) \qquad& \text{in } \mathbb{R}^n \times (0, T], \\ 
h(x,0) = 0 \qquad & \text{on } \mathbb{R}^n .
\end{cases}
$$
For $0<\ep\le1$ and $\sigma\in[0,1)$, there exists
$C=C(n,\sigma,T)>0$ such that
\[
|Dh|_{0,\sigma}
\le C\ep^{-\frac{1+\sigma}{2}}|\Gamma|_0.
\]
\end{lem}
\begin{proof}
The proof of the lemma is standard.
Let $\Phi_\varepsilon(x,t)=(4\pi \varepsilon t)^{-n/2} e^{-{|x|^2}/({4\varepsilon t})}$ be the fundamental solution of the heat equation. Because $h$ is a mild solution, 
Duhamel's principle gives 
\[
h(\cdot, t) = \int_0^t \Phi_\varepsilon(\cdot, t-s) *_x \Gamma(\cdot, s) \, ds.
\]
Since $\|\nabla \Phi_\varepsilon(\cdot, s)\|_{L^1} \leq C (\varepsilon s)^{-1/2}$ and $\|D^2_x \Phi_\varepsilon(\cdot, s)\|_{L^1} \leq C (\varepsilon s)^{-1}$, for any $0<\sigma<1$, we obtain
\begin{align*}
&\quad\,\,\|\nabla \Phi_\varepsilon(x_1 - \cdot, s) - \nabla \Phi_\varepsilon(x_2 - \cdot, s)\|_{L^1}\\
&\leq (2\|\nabla \Phi_\varepsilon( \cdot, s) \|_{L^1})^{1-\sigma} (\|D^2 \Phi_\varepsilon( \cdot, s) \|_{L^1}|x_1-x_2|)^\sigma\\
&\leq C (\varepsilon s)^{-\frac{1+\sigma}{2}} |x_1 - x_2|^\sigma.    
\end{align*}
Integrating this against $\|\Gamma(\cdot, s)\|_{L^\infty}$ over $s \in [0, t]$ gives 
\[
[D h(\cdot, t)]_{C^\sigma} \leq C \|\Gamma\|_{L^\infty} \varepsilon^{-\frac{1+\sigma}{2}} \int_0^t (t-s)^{-\frac{1+\sigma}{2}} \, ds\leq C_\sigma \|\Gamma\|_{L^\infty} \varepsilon^{-\frac{1+\sigma}{2}}  t^{\frac{1-\sigma}{2}}.
\]
The gradient kernel bound gives
\[
\|Dh(\cdot,t)\|_{L^\infty}
\le C\ep^{-1/2}t^{1/2}
\|\Gamma\|_{L^\infty(\R^n\times[0,t])}.
\]
This proves the case $\sigma=0$.
Taking the supremum over $t\in[0,T]$ and using
$\ep^{-1/2}\le\ep^{-(1+\sigma)/2}$ for $0<\ep\le1$
proves the full-norm estimate.
\end{proof}

Applying the preceding heat-kernel estimate to \eqref{eq:NMstep-new}
gives the following corollary.

\begin{cor}\label{C.5.7}
Let $0<\ep\le1$, $0<\sigma<1$, $\theta_k\ge1$, and let
$l\ge m\ge1$ be integers. 
If $h^k$ is a bounded classical solution of \eqref{eq:NMstep-new}
and the norms on the right-hand side are finite, the following
estimate holds, with an implicit constant depending only on
$n,l,m,\sigma,T$, and $\rho$:
\begin{align*}
 |h^k|_{l+1,\sigma}\lesssim \ep^{-\frac{1+\sigma}{2}}\left(\theta_k^{l-m}|R_k|_m+|D_pH(x,D\widetilde u^k)|_{l} |h^k|_1 + |D_pH(x,D\widetilde u^k)|_0 |h^k|_{l+1} \right).     
\end{align*}
\end{cor}
\begin{proof}
Writing \eqref{eq:NMstep-new} as a heat equation and differentiating
its Duhamel formula in space through order $l$, Lemma~\ref{L.7.6} gives
\begin{align*}
|Dh^k|_{l,\sigma}&\lesssim \ep^{-\frac{1+\sigma}{2}}(|S_{\theta_k}R_k|_l +|D_pH(x,D\widetilde u^k)\cdot Dh^k|_l).    
\end{align*}
By \eqref{eq:smoothing2-new} and product estimates,
\[
|Dh^k|_{l,\sigma}\lesssim \ep^{-\frac{1+\sigma}{2}}\left(\theta_k^{l-m}|R_k|_m+|D_pH(x,D\widetilde u^k)|_{l} |h^k|_1 + |D_pH(x,D\widetilde u^k)|_0 |h^k|_{l+1} \right).
\]
Finally, the comparison principle gives
\[
|h^k|_0
\le T|S_{\theta_k}R_k|_0
\le T\|\rho\|_{L^1}|R_k|_m.
\]
Since $|h^k|_{l+1,\sigma}\lesssim|h^k|_0+|Dh^k|_{l,\sigma}$,
$\theta_k\ge1$, and $\ep\le1$, this completes the proof.
\end{proof}

\subsection{The one-step estimate}

The following identity separates the new residual into the
residual-smoothing error, the drift-smoothing error, and the
quadratic Newton remainder.

\begin{lem}\label{lem:residual-identity-new}
Let $u^+=u+h$, $R=f_\ep[u]$, and
\[
    \widetilde u=\bar g+S_\theta(u-\bar g).
\]
Assume that $h$ solves
\begin{equation}\label{eq:abstract-step-new}
\begin{cases}
    \Lop_{\widetilde u}^\ep h=-S_\theta R,\\
    h(x,0)=0.
\end{cases}
\end{equation}
Then
\begin{equation}\label{eq:residual-identity-new}
    f_\ep[u^+]
    =(I-S_\theta)R+\bigl[D_pH(x,Du)-D_pH(x,D\widetilde u)\bigr]\cdot Dh+\mathcal Q_u(h),
\end{equation}
where
\begin{equation*}%\label{eq:Qdef-new}
    \mathcal Q_u(h):=H(x,Du+Dh)-H(x,Du)-D_pH(x,Du)\cdot Dh.
\end{equation*}
Equivalently,
\begin{equation}\label{eq:Qintegral-new}
    \mathcal Q_u(h)=\int_0^1(1-r)D^2_{pp}H(x,Du+rDh)[Dh,Dh]\,dr.
\end{equation}
\end{lem}

\begin{proof}
By Taylor's expansion,
\begin{align*}
    f_\ep[u+h]
    &=u_t+h_t+H(x,Du+Dh)-\ep\Delta u-\ep\Delta h\\
    &=R+h_t+D_pH(x,Du)\cdot Dh-\ep\Delta h+\mathcal Q_u(h).
\end{align*}
Using
\[
    h_t+D_pH(x,D\widetilde u)\cdot Dh-\ep\Delta h=-S_\theta R,
\]
we obtain \eqref{eq:residual-identity-new}.
\end{proof}

\begin{prop}[One-step Nash--Moser estimate]\label{prop:onestep-new}
Let $m\ge1$ and $N\ge m+2$ be integers, let $\sigma\in[0,1)$,
and let $0<\ep\le1$ and $\theta\ge1$.
Choose the convolution kernel above with
$N+\sigma\le N_*+1$.
Assume that $H\in {C_b^{m+2}}(\R^n\times B(0,{r}))$ for every ${r}>0$ and that $g\in C_b^{N+\sigma}(\R^n)$. Suppose that $u\in C^{2,1}(Q_T)\cap C_b(\overline Q_T)$ satisfies
\begin{equation}\label{eq:low-u-bound-new}
    |u|_{m+1}\le M\quad\text{and}\quad |u|_{N,\sigma}<\infty,
\end{equation}
and assume that $R=f_\ep[u]\in C_b(\overline Q_T)$ satisfies
$|R|_1\le1$ and $|R|_{N,\sigma}<\infty$.
Let $\widetilde u=\bar g+S_\theta(u-\bar g)$, 
and let $h\in C^{2,1}(Q_T)\cap C_b(\overline Q_T)$ solve
\eqref{eq:abstract-step-new}.
Then
\begin{equation}\label{eq:h-estimate-new}
    |h|_{m+1}\le C\theta |R|_m,
\end{equation}
and
\begin{equation}\label{eq:onestep-estimate-new}
    |f_\ep[u+h]|_m
    \le C\theta^{m-N-\sigma}|R|_{N,\sigma}
        +C\theta^{m+1-N-\sigma}|u-\bar g|_{N,\sigma}|R|_m
        +C\theta|R|_m^2.
\end{equation}
Here $C$ depends on $n,m,N,\sigma,M,T,\rho,|g|_{m+2}$, and
suitable local bounds for derivatives of $H$ through order $m+2$,
but is independent of $0<\ep\le1$, $\theta\ge1$, and
$|u|_{N,\sigma}$.
\end{prop}

\begin{proof}
Set
\[
    b=D_pH(x,D\widetilde u).
\]
Because $|u|_{m+1}\le M$, the smoothing estimates give $|\widetilde u|_{m+1}\le C(M,g)$ and hence, also using that $m\geq 1$,
\begin{equation}\label{b1}
|b|_1\le C(M,g,H).
\end{equation}
Moreover, by \eqref{eq:smoothing2-new} and \eqref{eq:low-u-bound-new},
\[
|b|_{m+1}
\le C\bigl(1+|\widetilde u|_{m+2}\bigr)
\le C\bigl(1+|g|_{m+2}+\theta|u-\bar g|_{m+1}\bigr)
\le C\theta.
\]

Applying Lemma \ref{lem:linear-new} with $s=m+1$ and $f=-S_\theta R$, and using
\[
    |S_\theta R|_{m+1}\le C\theta |R|_m,
    \qquad |S_\theta R|_1\le C|R|_1\leq C|R|_m,
\]
we obtain
\[
    |h|_{m+1}\le  C\theta |R|_m.
\]
This proves \eqref{eq:h-estimate-new}.  Similarly, by Lemma \ref{lem:linear-new}, 
\begin{equation}\label{eq:h1-estimate-new}
    |h|_1\le C|R|_1\le C|R|_m.
\end{equation}

We now use the residual identity \eqref{eq:residual-identity-new}.  The first term on the right-hand side of \eqref{eq:residual-identity-new} is controlled by \eqref{eq:smoothing3-new}:
\begin{equation}\label{5101}
    |(I-S_\theta)R|_m\le C\theta^{m-N-\sigma}|R|_{N,\sigma}.    
\end{equation}
For the drift-smoothing error, the mean-value formula and product estimates give
\begin{align*}
&\left|\bigl[D_pH(x,Du)-D_pH(x,D\widetilde u)\bigr]\cdot Dh\right|_m \\
&\qquad\le C|u-\widetilde u|_{m+1}|h|_1+C|u-\widetilde u|_1|h|_{m+1}.
\end{align*}
Notice that
\[
u-\widetilde u= (I-S_\theta)(u-\bar g).
\]
Using \eqref{eq:smoothing3-new}, \eqref{eq:h-estimate-new},  and \eqref{eq:h1-estimate-new},
\begin{align*}
    &|u-\widetilde u|_{m+1}\le C\theta^{m+1-N-\sigma}|u-\bar g|_{N,\sigma},\quad
    |u-\widetilde u|_1\le C\theta^{1-N-\sigma}|u-\bar g|_{N,\sigma},\\
    &|h|_{m+1}\le C\theta|R|_m, \qquad
    |h|_1\le C|R|_m.
\end{align*}
Thus, since $m\ge1$ and $\theta\ge1$, we obtain
\begin{equation}\label{5102}
    \left|\bigl[D_pH(x,Du)-D_pH(x,D\widetilde u)\bigr]\cdot Dh\right|_m
    \le C\theta^{m+1-N-\sigma}|u-\bar g|_{N,\sigma}|R|_m.
\end{equation}

Finally, \eqref{eq:h1-estimate-new} and $|R|_1\le1$ give
$|Du+rDh|_0\le C$ uniformly for $r\in[0,1]$.
The composition estimate \eqref{eq:tame-composition} therefore gives
\[
|D^2_{pp}H(x,Du+rDh)|_m
\le C\bigl(1+|u|_{m+1}+|h|_{m+1}\bigr)
\le C(1+\theta|R|_m),
\]
and
\[
|D^2_{pp}H(x,Du+rDh)|_0\le C.
\]
Hence, by \eqref{eq:Qintegral-new}, product estimates,  and \eqref{eq:h1-estimate-new}, the last term of \eqref{eq:residual-identity-new} is bounded by
\begin{equation}\label{5103}
\begin{aligned}
 |\mathcal Q_u(h)|_m &\le C\sup_{r\in[0,1]}(|D^2_{pp}H(x,Du+rDh)|_m|h|_1^2+|D^2_{pp}H(x,Du+rDh)|_0|h|_{m+1}|h|_1)\\
&\le C(1+\theta|R|_m)|R|_1^2+C\theta|R|_m|R|_1\leq C\theta|R|_m^2,
\end{aligned}
\end{equation}
where the last inequality is due to $|R|_1\leq 1$, $|R|_1\leq |R|_m$ and $\theta\geq 1$.
Combining the three estimates of \eqref{5101}, \eqref{5102}, and \eqref{5103} gives \eqref{eq:onestep-estimate-new}.
\end{proof}

\subsection{Proof of the convergence criterion}

%{\color{teal}

\begin{proof}[Proof of Theorem \ref{thm:criterion-new}]
Set $N:=m+3$ and $a:=3/2$, and let $\sigma\in(0,1)$, $p:=(1+\sigma)/2$,
and $\chi\in(1,4/3)$ be parameters depending only on $\mu$, to be
specified in Step~3.

First, for each fixed $\ep>0$ and $\theta_0\ge2$, the iteration is well
defined and the high norms used below are finite. Indeed, the assumptions
imply $u^0\in X_T$ and $g\in C_b^{2+\alpha}(\R^n)$. If $u^k\in X_T$,
then $\widetilde u^k\in X_T$ and
\[
    D_pH(x,D\widetilde u^k),\ S_{\theta_k}R_k
    \in C_b^{\alpha,\alpha/2}(\overline Q_T).
\]
Here the bounds on $H$ apply on the bounded gradient range of
$\widetilde u^k$. Parabolic Schauder theory therefore gives a unique
solution $h^k\in X_T$ of \eqref{eq:NMstep-new}. Thus $u^{k+1}\in X_T$,
and induction proves that the sequence is well defined. Moreover, the
regularity assumptions give
\[
    |u^0|_{N+1,\sigma}
    +\sup_{0<\ep\le1}|f_\ep[u^0]|_{N,\sigma}<\infty.
\]
The finiteness of $|u^k|_{N+1,\sigma}$ and $|R_k|_{N,\sigma}$ then
propagates by induction: Lemma~\ref{lem:linear-new} with $s=N+1$ first
gives $|h^k|_{N+1}<\infty$, Corollary~\ref{C.5.7} with $l=N$ then gives
$|h^k|_{N+1,\sigma}<\infty$, and the update and residual identity give
the corresponding bounds at the next step. No estimate uniform in $\ep$
is needed for these qualitative properties. We divide the rest of the
proof into four steps.

\medskip

\noindent {\bf Step 1.} The inductive bootstrap and low-norm bound.

We now prove quantitative residual decay together with a uniform
low-norm bound. Let
\[
    M:=1+|u^0|_{m+1}.
\]
We will prove, for $\Theta_*$ sufficiently large (chosen after $\sigma$
and $\chi$ are fixed in Step~3) and every admissible
$\theta_0\ge\Theta_*$, that for all $k\ge0$,
\begin{equation}\label{eq:NM-induction-low}
    |u^k|_{m+1}\le M,
\end{equation}
and
\begin{equation}\label{eq:NM-induction-residual}
    |R_k|_m\le \theta_k^{-a}.
\end{equation}
The case $k=0$ follows from the choice of $M$ and from
\eqref{eq:initial-small-new}. Assume that \eqref{eq:NM-induction-low} and
\eqref{eq:NM-induction-residual} hold for $0,\ldots,k$.

By Proposition~\ref{prop:onestep-new}, applied with the constant
corresponding to the bound $M$, we have
\begin{equation}\label{eq:h-low-corrected}
    |h^i|_{m+1}
    \le C\theta_i |R_i|_m
    \le C\theta_i^{1-a},
    \qquad 0\le i\le k.
\end{equation}
Since $a>1$ and $\chi>1$, taking $\Theta_*$ sufficiently large ensures
$\sum_{i=0}^{\infty} C\theta_i^{1-a}\le 1$. Therefore
\[
    |u^{k+1}|_{m+1}
    \le |u^0|_{m+1}+\sum_{i=0}^{k}|h^i|_{m+1}
    \le M.
\]
Thus the low-norm bootstrap closes.

\medskip

\noindent {\bf Step 2.} High-norm estimates of $h^i$, $u^i$, and $R_i$.

To close the residual induction, the one-step estimate requires
bounds on $|u^k-\bar g|_{N,\sigma}$ and $|R_k|_{N,\sigma}$.
We derive these from the correction estimates below.
These high norms may grow; their growth must be outweighed
by the smoothing factors.
Set
\[
    b^i:=D_pH(x,D\widetilde u^i).
\]
Since
\[
    \widetilde u^i=\bar g+S_{\theta_i}(u^i-\bar g),
\]
the smoothing estimates \eqref{eq:smoothing1-new}--\eqref{eq:smoothing2-new}
and the local $C_b^{N+2}$ bounds on $H$ imply, for every integer
$0\le s\le N+1$,
\begin{equation*}
    |b^i|_s\le C\theta_i^{(s-m)_+}.
\end{equation*}
Here $x_+=\max\{x,0\}$, and the constants depend on $M,g,H,s,m$, but not
on $\ep$, $\theta_i$, or $i$.

We also use Lemma~\ref{lem:linear-new}. Since
\[
\begin{cases}
    h^i_t+b^i\cdot Dh^i-\ep\Delta h^i=-S_{\theta_i}R_i,\\
    h^i(\cdot,0)=0,
\end{cases}
\]
Lemma~\ref{lem:linear-new} gives, for integers $m\le s\le N+1$,
\begin{equation*}
\begin{aligned}
    |h^i|_s
    &\le C\left(
        |S_{\theta_i}R_i|_s
        +|b^i|_s |S_{\theta_i}R_i|_1
    \right) \\
    &\le C\theta_i^{s-m}|R_i|_m
    \le C\theta_i^{s-m-a}.
\end{aligned}
\end{equation*}
At first order, the same lemma gives
\begin{equation}\label{eq:h-one-corrected}
    |h^i|_1\le C\theta_i^{-a}.
\end{equation}

We now use Corollary~\ref{C.5.7}. For every integer $m+1\le\ell\le N+1$,
using the above estimates,
\[
\begin{aligned}
    |h^i|_{\ell,\sigma}
    &\le C\ep^{-p}
    \left(
        \theta_i^{\ell-1-m}|R_i|_m
        +|b^i|_{\ell-1}|h^i|_1
        +|b^i|_0|h^i|_\ell
    \right) \\
    &\le C\ep^{-p}
    \left(
        \theta_i^{\ell-1-m-a}
        +\theta_i^{\ell-1-m-a}
        +\theta_i^{\ell-m-a}
    \right) \\
    &\le C\ep^{-p}\theta_i^{\ell-m-a}.
\end{aligned}
\]
In particular, since $N=m+3$,
\begin{equation}\label{eq:h-N-sigma-corrected}
    |h^i|_{N,\sigma}
    \le C\ep^{-p}\theta_i^{3-a},
\end{equation}
and
\begin{equation}\label{eq:h-N1-sigma-corrected}
    |h^i|_{N+1,\sigma}
    \le C\ep^{-p}\theta_i^{4-a}.
\end{equation}
Summing over $i=0,\ldots,k-1$ and using $a<2$ and
$|u^0|_{N+1,\sigma}<\infty$ yields
\begin{equation}\label{eq:u-N-sigma-corrected}
    |u^k-\bar g|_{N,\sigma}
    \le
    C\left(
        1+\ep^{-p}\theta_k^{(3-a)/\chi}
    \right),
\end{equation}
and
\begin{equation}\label{eq:u-N1-sigma-corrected}
    |u^k|_{N+1,\sigma}
    \le
    C\left(
        1+\ep^{-p}\theta_k^{(4-a)/\chi}
    \right).
\end{equation}

We now estimate the high residual directly from the residual identity.
For each $i$, define
\[
    A^i:=
    D_pH(x,Du^i)-D_pH(x,D\widetilde u^i).
\]
By the residual identity \eqref{eq:residual-identity-new},
\begin{equation}\label{eq:residual-identity-used-corrected}
    R_{i+1}
    =
    (I-S_{\theta_i})R_i
    +
    A^i\cdot Dh^i
    +
    \mathcal Q_{u^i}(h^i).
\end{equation}
We estimate the terms in $C^{N,\sigma}$. First, by \eqref{eq:smoothing1-new},
\[
    |(I-S_{\theta_i})R_i|_{N,\sigma}
    \le C_\rho|R_i|_{N,\sigma},
    \qquad C_\rho:=1+\|\rho\|_{L^1}.
\]
Next, since $u^i-\widetilde u^i=(I-S_{\theta_i})(u^i-\bar g)$ and
$|u^i-\bar g|_{m+1}\le C$, the smoothing estimate \eqref{eq:smoothing3-new}
gives
\begin{equation}\label{eq:A-zero-corrected}
    |A^i|_0
    \le C|u^i-\widetilde u^i|_1
    \le C\theta_i^{-m}|u^i-\bar g|_{m+1}
    \le C\theta_i^{-m}.
\end{equation}
Moreover, by the composition estimates, the local $C_b^{N+2}$
bounds on $H$, and \eqref{eq:u-N1-sigma-corrected},
\begin{equation}\label{eq:A-high-corrected}
    |A^i|_{N,\sigma}
    \le
    C\left(
        1+\ep^{-p}\theta_i^{(4-a)/\chi}
    \right).
\end{equation}
Therefore, using \eqref{eq:h-one-corrected}, \eqref{eq:h-N1-sigma-corrected},
\eqref{eq:A-zero-corrected}, and \eqref{eq:A-high-corrected},
\[
\begin{aligned}
    |A^i\cdot Dh^i|_{N,\sigma}
    &\le C\left(
        |A^i|_{N,\sigma}|h^i|_1
        +|A^i|_0|h^i|_{N+1,\sigma}
    \right) \\
    &\le C\left(
        \theta_i^{-a}
        +\ep^{-p}\theta_i^{(4-a)/\chi-a}
        +\ep^{-p}\theta_i^{4-a-m}
    \right).
\end{aligned}
\]
Writing $a\wedge m:=\min\{a,m\}$ and using $\theta_i\ge1$, $\ep\le1$, and
$\chi>1$, this implies
\begin{equation}\label{eq:A-Dh-corrected}
    |A^i\cdot Dh^i|_{N,\sigma}
    \le
    C\left(
        1+\ep^{-p}\theta_i^{4-a-a\wedge m}
    \right).
\end{equation}

For the quadratic term, applying the product and composition
estimates to \eqref{eq:Qintegral-new}, using the local $C_b^{N+3}$
bounds on $H$, gives
\[
\begin{aligned}
    |\mathcal Q_{u^i}(h^i)|_{N,\sigma}
    &\le
    C\left(
        1+|u^i|_{N+1,\sigma}+|h^i|_{N+1,\sigma}
    \right)|h^i|_1^2
    +
    C|h^i|_{N+1,\sigma}|h^i|_1.
\end{aligned}
\]
Using \eqref{eq:h-one-corrected}, \eqref{eq:h-N1-sigma-corrected}, and
\eqref{eq:u-N1-sigma-corrected}, we obtain
\[
\begin{aligned}
    |\mathcal Q_{u^i}(h^i)|_{N,\sigma}
    &\le
    C\left(
        \theta_i^{-2a}
        +\ep^{-p}\theta_i^{(4-a)/\chi-2a}
        +\ep^{-p}\theta_i^{4-2a}
    \right) \\
    &\le
    C\left(
        1+\ep^{-p}\theta_i^{4-2a}
    \right).
\end{aligned}
\]
Combining this with \eqref{eq:residual-identity-used-corrected} and
\eqref{eq:A-Dh-corrected} gives
\[
    |R_{i+1}|_{N,\sigma}
    \le
    C_\rho|R_i|_{N,\sigma}
    +
    C\left(
        1+\ep^{-p}\theta_i^{4-a-a\wedge m}
    \right).
\]
Set $d:=4-a-(a\wedge m)>0$. Since $\ep^{-p}\theta_i^d\ge1$, iterating
the recurrence gives
\[
    |R_k|_{N,\sigma}
    \le C_\rho^k|R_0|_{N,\sigma}
    +C\ep^{-p}\sum_{i=0}^{k-1}C_\rho^{k-1-i}\theta_i^d.
\]
For $k\ge1$, writing $j=k-1-i$ and using $\theta_0\ge2$, we obtain
\[
\begin{aligned}
    \sum_{i=0}^{k-1}C_\rho^{k-1-i}\theta_i^d
    &\le
    \theta_{k-1}^d
    \sum_{j=0}^{\infty}C_\rho^j\,2^{-d(\chi^j-1)}\\
    &\le C\theta_k^{d/\chi}.
\end{aligned}
\]
The series converges because $\chi>1$. Moreover, the uniform initial
bound gives, for every $\eta>0$,
\[
    C_\rho^k|R_0|_{N,\sigma}\le C_\eta\theta_k^\eta,
\]
with constants independent of $\ep,\theta_0,k$. Together with the
initial bound when $k=0$, this yields
\begin{equation}\label{eq:R-high-corrected}
    |R_k|_{N,\sigma}
    \le
    C_\eta \theta_k^\eta
    +C\ep^{-p}\theta_k^{(4-a-a\wedge m)/\chi}.
\end{equation}

\medskip

\noindent {\bf Step 3.} Closing the residual induction and parameter selection.

We now apply Proposition~\ref{prop:onestep-new}. Since $N=m+3$,
\eqref{eq:onestep-estimate-new}, \eqref{eq:R-high-corrected},
\eqref{eq:u-N-sigma-corrected}, and the induction hypothesis give
\[
\begin{aligned}
    |R_{k+1}|_m
    &\le
    C\theta_k^{-3-\sigma}|R_k|_{N,\sigma}
    +
    C\theta_k^{-2-\sigma}|u^k-\bar g|_{N,\sigma}|R_k|_m
    +
    C\theta_k|R_k|_m^2 \\
    &\le
    C_\eta\theta_k^{-3-\sigma+\eta}
    +
    C\ep^{-p}
        \theta_k^{-3-\sigma+(4-a-a\wedge m)/\chi} \\
    &\quad
    +
    C\theta_k^{-2-\sigma-a}
    +
    C\ep^{-p}\theta_k^{-2-\sigma-a+(3-a)/\chi}
    +
    C\theta_k^{1-2a}.
\end{aligned}
\]
Using $\ep\ge\theta_0^{-\mu}\ge\theta_k^{-\mu}$, this becomes
\begin{equation}\label{eq:R-k1-five-terms}
\begin{aligned}
    |R_{k+1}|_m
    &\le
    C_\eta\theta_k^{-3-\sigma+\eta}
    +
    C\theta_k^{-3-\sigma+p\mu+(4-a-a\wedge m)/\chi} \\
    &\quad
    +
    C\theta_k^{-2-\sigma-a}
    +
    C\theta_k^{-2-\sigma-a+p\mu+(3-a)/\chi}
    +
    C\theta_k^{1-2a}.
\end{aligned}
\end{equation}

We now specify $\sigma$, $\eta$, and $\chi$, depending only on $\mu$, so
that every exponent in \eqref{eq:R-k1-five-terms} is strictly less than
$-a\chi$. 
For the second term in \eqref{eq:R-k1-five-terms}, the worst case
is \(m=1\). As \(\chi\downarrow1\), its exponent tends to
\(-a-\sigma+p\mu\), whereas the target exponent tends to \(-a\).
We therefore choose, depending only on \(\mu\),
\[
    \sigma\in\left(\frac{\mu}{2-\mu},1\right),
    \qquad
    p=\frac{1+\sigma}{2},
\]
which guarantees that $\sigma-p\mu>0$. Next, fix $0<\eta<\sigma$,
depending only on $\mu$. With $a=3/2$, we have
\[
    3+\sigma-p\mu-\frac{3-a}{\chi}-a\chi
    \longrightarrow \sigma-p\mu>0
    \qquad\text{as }\chi\downarrow1.
\]
We may therefore choose $\chi\in(1,4/3)$, depending only on $\mu$,
sufficiently close to $1$ so that
\[
    3+\sigma-p\mu-\frac{3-a}{\chi}>a\chi.
\]
Since $\chi<4/3$, $a=3/2$, and $a\wedge m\ge1$ for every integer $m\ge1$,
this choice implies that for every $m\ge1$,
\begin{equation}\label{eq:param-choice-corrected}
    a-1>a(\chi-1),
    \qquad
    3+\sigma-\eta>a\chi,
\end{equation}
\begin{equation}\label{eq:param-choice-corrected-2}
    3+\sigma-p\mu-\frac{4-a-(a\wedge m)}{\chi}
    \ge
    3+\sigma-p\mu-\frac{3-a}{\chi}
    >a\chi,
\end{equation}
and
\begin{equation}\label{eq:param-choice-corrected-3}
    2+\sigma+a-p\mu-\frac{3-a}{\chi}
    =
    \left(3+\sigma-p\mu-\frac{3-a}{\chi}\right)+(a-1)
    >a\chi.
\end{equation}
By \eqref{eq:param-choice-corrected}--\eqref{eq:param-choice-corrected-3},
there exists $\omega>0$ such that each exponent on the right-hand side
of \eqref{eq:R-k1-five-terms} is at most $-a\chi-\omega$. Hence
\[
    |R_{k+1}|_m
    \le
    C_\eta \theta_k^{-a\chi-\omega}.
\]
Choose $\Theta_*$ large enough to satisfy the requirement in Step~1 and
$C_\eta\Theta_*^{-\omega}\le1$. Then every admissible $\theta_0\ge\Theta_*$
satisfies
\[
    |R_{k+1}|_m
    \le
    \theta_k^{-a\chi}
    =
    \theta_{k+1}^{-a}.
\]
This closes the induction.

\medskip

\noindent {\bf Step 4.} Convergence.

It remains to prove convergence. From \eqref{eq:h-low-corrected},
\[
    |u^{k+1}-u^k|_{m+1}
    =
    |h^k|_{m+1}
    \le
    C\theta_k^{1-a}.
\]
Since $a>1$,
\[
    \sum_{k=0}^\infty \theta_k^{1-a}
    =
    \sum_{k=0}^\infty
    \left(\theta_0^{1-a}\right)^{\chi^k}
    <\infty.
\]
Thus $u^k$ converges with respect to $|\cdot|_{m+1}$ to some function $u$
satisfying $|u|_{m+1}\le M$, and $u(x,0)=g(x)$ because each $h^k$ has zero
initial trace.

Since $|R_k|_m\to0$ and $|u^k-u|_{m+1}\to0$, we have
\[
    |H(x,Du^k)-H(x,Du)|_{m-1}\to0,
    \qquad
    |\ep\Delta u^k-\ep\Delta u|_{m-1}\to0.
\]
Moreover,
\[
    (u^k)_t
    =
    R_k-H(x,Du^k)+\ep\Delta u^k
\]
converges with respect to $|\cdot|_{m-1}$, and in particular uniformly on
$\overline Q_T$. Since $u^k\to u$ uniformly on $\overline Q_T$, the limit
$u$ is continuously differentiable in time and satisfies
\[
    u_t+H(x,Du)-\ep\Delta u=0
    \quad\text{in }Q_T,
    \qquad
    u(x,0)=g(x)
    \quad\text{on }\R^n.
\]
Since $m\ge1$, $u$ is therefore a bounded classical solution of
\eqref{eq:viscousHJ-eps}.

Finally,
\[
\begin{aligned}
    |u^k-u|_{m+1}
    &\le
    \sum_{j\ge k}|u^{j+1}-u^j|_{m+1} \\
    &\le
    C\sum_{j\ge k}\theta_j^{1-a}
    \le
    C\theta_k^{1-a}
    =
    C\left(\theta_0^{1-a}\right)^{\chi^k}.
\end{aligned}
\]
This completes the proof.
\end{proof}

\section{PINN residual certificates for inexact Newton steps}
\label{sec:pinn-certificate}

In the preceding sections, we analyzed the exact Newton iteration where each step requires solving a linear uniformly parabolic equation for the Newton correction $h^k := u^{k+1} - u^k$. 
These linear subproblems generally require numerical approximation. In high-dimensional settings, 
a natural approach is to parameterize the correction $h^k$ at each step by a physics-informed neural network (PINN). 
Optimization and representation errors lead to an inexact Newton scheme.

We bound the error of the updated iterate using the continuous
linearized residual and a bound on the quadratic Newton remainder,
without using the exact solution.
For simplicity, we work on the periodic cylinder
\[
    Q_T:=\T^n\times(0,T).
\]
An analogous argument applies on bounded smooth domains, provided the boundary
condition is imposed exactly or the boundary mismatch and corresponding error estimates are taken into account.  We consider
\begin{equation}\label{eq:pinn-hj}
\begin{cases}
    u_t+H(x,Du)-\Delta u=0
   \qquad &\text{in } Q_T,\\
    u(x,0)=g(x)
    \qquad&\text{on } \T^n .
\end{cases}
\end{equation}
Set
\begin{equation}\label{eq:pinn-f-def}
    f[v]:=v_t+H(x,Dv)-\Delta v .
\end{equation}
Throughout this section, we assume that the Hamiltonian has been modified
outside the relevant gradient range so that
\eqref{eq:conditions-H-modified} holds. Set
\begin{equation}\label{eq:pinn-bounded-H}
    M_1:=\norm{D_pH}_{L^\infty(\T^n\times\R^n)}<\infty,
    \qquad
    M_2:=\sup_{(x,p)\in\T^n\times\R^n}
    \|D^2_{pp}H(x,p)\|_{\mathrm{op}}<\infty,
\end{equation}
where \(\|\cdot\|_{\mathrm{op}}\) denotes the Euclidean operator norm.

We use the parabolic energy norm
\begin{equation}\label{eq:pinn-energy-norm}
    \norm{z}_{\cE(Q_T)}^2
    :=
    \sup_{0\le t\le T}
    \norm{z(\cdot,t)}_{L^2(\T^n)}^2
    +
    \int_0^T
    \norm{Dz(\cdot,t)}_{L^2(\T^n)}^2\,dt .
\end{equation}

\subsection{Inexact Newton correction}

Let \(\widehat u^k\) be the current approximation and assume that it satisfies
the initial condition exactly:
\[
    \widehat u^k(x,0)=g(x).
\]
Define the nonlinear residual and frozen Newton drift by
\begin{equation}\label{eq:pinn-current-residual}
    \widehat R_k:=f[\widehat u^k],
    \qquad
    \widehat b^k:=D_pH(x,D\widehat u^k).
\end{equation}
The exact Newton correction associated with \(\widehat u^k\) solves
\begin{equation}\label{eq:pinn-exact-correction}
\begin{cases}
    (h_*^k)_t+\widehat b^k\cdot Dh_*^k-\Delta h_*^k
    =
    -\widehat R_k
    \qquad&\text{in } Q_T,\\
    h_*^k(x,0)=0
    \qquad&\text{on } \T^n .
\end{cases}
\end{equation}
In the PINN implementation, we replace \(h_*^k\) by a neural correction
\(h^k\).  One may impose the zero initial trace exactly by using the ansatz
\[
    h^k(x,t)=t\,\mathcal N^k(x,t),
\]
where \(\mathcal N^k\) is a smooth neural network, periodic in \(x\).  
The trained correction is measured by its continuous linearized residual
\begin{equation}\label{eq:pinn-rho-natural}
    \rho_k
    :=
    (h^k)_t
    +
    \widehat b^k\cdot Dh^k
    -
    \Delta h^k
    +
    \widehat R_k .
\end{equation}
The updated iterate is
\begin{equation}\label{eq:pinn-update}
    \widehat u^{k+1}:=\widehat u^k+h^k .
\end{equation}
If \(\rho_k=0\), then \(h^k\) is the exact Newton correction for the frozen
iterate \(\widehat u^k\); otherwise Theorem~\ref{thm:pinn-aposteriori} below
quantifies the effect of the inexact solve.

\subsection{A posteriori certificate}

We now give the deterministic certificate.  The proof uses three elementary
ingredients: a linear parabolic energy estimate, nonlinear residual stability,
and the residual identity for one inexact Newton update.

\begin{lem}
\label{lem:pinn-linear-energy}
Let \(b\in L^\infty(Q_T;\R^n)\) and \(\eta\in L^2(Q_T)\), and let \(z\)
be a sufficiently smooth periodic solution of
\begin{equation}\label{eq:pinn-linear-model}
\begin{cases}
    z_t+b\cdot Dz-\Delta z=\eta
    \qquad&\text{in } Q_T,\\
    z(x,0)=0
    \qquad&\text{on } \T^n .
\end{cases}
\end{equation}
Then
\begin{equation}\label{eq:pinn-linear-energy-estimate}
    \norm{z}_{\cE(Q_T)}
    \le
    C\bigl(T,\norm{b}_{L^\infty(Q_T)}\bigr)
    \norm{\eta}_{L^2(Q_T)} .
\end{equation}
\end{lem}

\begin{proof}
Set \(B:=\norm{b}_{L^\infty(Q_T)}\).  Multiplying
\eqref{eq:pinn-linear-model} by \(z\), integrating over \(\T^n\), and using
periodicity, we obtain
\[
    \frac12\frac{d}{dt}\norm{z(\cdot,t)}_{L^2(\T^n)}^2
    +
    \norm{Dz(\cdot,t)}_{L^2(\T^n)}^2
    =
    \int_{\T^n}\eta z\,dx
    -
    \int_{\T^n} b\cdot Dz\,z\,dx .
\]
By Young's inequality,
\[
    \int_{\T^n}\eta z\,dx
    \le
    \frac12\norm{\eta(\cdot,t)}_{L^2(\T^n)}^2
    +
    \frac12\norm{z(\cdot,t)}_{L^2(\T^n)}^2,
\]
and
\[
    \left|
    \int_{\T^n} b\cdot Dz\,z\,dx
    \right|
    \le
    \frac12\norm{Dz(\cdot,t)}_{L^2(\T^n)}^2
    +
    \frac{B^2}{2}\norm{z(\cdot,t)}_{L^2(\T^n)}^2 .
\]
Hence
\[
    \frac{d}{dt}\norm{z(\cdot,t)}_{L^2(\T^n)}^2
    +
    \norm{Dz(\cdot,t)}_{L^2(\T^n)}^2
    \le
    \norm{\eta(\cdot,t)}_{L^2(\T^n)}^2
    +
    (1+B^2)\norm{z(\cdot,t)}_{L^2(\T^n)}^2 .
\]
Since \(z(\cdot,0)=0\), Gronwall's inequality gives
\[
    \sup_{0\le t\le T}
    \norm{z(\cdot,t)}_{L^2(\T^n)}^2
    \le
    e^{(1+B^2)T}
    \norm{\eta}_{L^2(Q_T)}^2 .
\]
Integrating the differential inequality in time gives the corresponding
gradient estimate.  Combining the two bounds proves
\eqref{eq:pinn-linear-energy-estimate}.
\end{proof}

\begin{lem}
\label{lem:pinn-nonlinear-stability}
Let \(u\) be the classical periodic solution of \eqref{eq:pinn-hj}.  Let
\(v\) be a sufficiently smooth periodic function satisfying
\[
    v(x,0)=g(x)
    \qquad\text{on } \T^n,
\]
and suppose
\begin{equation}\label{eq:pinn-v-residual}
    f[v]=\eta
    \qquad\text{in } Q_T .
\end{equation}
Assume \eqref{eq:pinn-bounded-H}.  Then
\begin{equation}\label{eq:pinn-nonlinear-stability}
    \norm{v-u}_{\cE(Q_T)}
    \le
    C(T,M_1)\norm{\eta}_{L^2(Q_T)} .
\end{equation}
\end{lem}

\begin{proof}
Set \(e:=v-u\).  Since \(v\) and \(u\) have the same initial data,
\(e(x,0)=0\).  Moreover, \(f[v]=\eta\) and \(f[u]=0\) imply
\[
    e_t+H(x,Dv)-H(x,Du)-\Delta e=\eta .
\]
We have
\[
    H(x,Dv)-H(x,Du)=a(x,t)\cdot De,
\]
where
\[
    a(x,t)
    :=
    \int_0^1D_pH\bigl(x,Du+s(Dv-Du)\bigr)\,ds .
\]
By \eqref{eq:pinn-bounded-H}, \(\norm{a}_{L^\infty(Q_T)}\le M_1\).  Thus
\(e\) solves
\[
\begin{cases}
    e_t+a\cdot De-\Delta e=\eta,\\
    e(x,0)=0. 
\end{cases}
\]
Lemma \ref{lem:pinn-linear-energy} gives \eqref{eq:pinn-nonlinear-stability}.
\end{proof}

\begin{lem}
\label{lem:pinn-residual-identity}
Let \(\widehat u^k\) and \(h^k\) be sufficiently smooth periodic functions
satisfying
\[
    \widehat u^k(x,0)=g(x),
    \qquad
    h^k(x,0)=0
    \qquad\text{on } \T^n .
\]
Set \(\widehat u^{k+1}:=\widehat u^k+h^k\), define
\[
    \widehat R_k:=f[\widehat u^k],
    \qquad
    \widehat b^k:=D_pH(x,D\widehat u^k),
\]
and let \(\rho_k\) be given by \eqref{eq:pinn-rho-natural}.  Then
\begin{equation}\label{eq:pinn-residual-identity}
    f[\widehat u^{k+1}]
    =
    \rho_k+\mathcal Q_k,
\end{equation}
where
\begin{equation}\label{eq:pinn-Q-natural}
    \mathcal Q_k
    :=
    H(x,D\widehat u^k+Dh^k)
    -
    H(x,D\widehat u^k)
    -
    D_pH(x,D\widehat u^k)\cdot Dh^k .
\end{equation}
Moreover, under \eqref{eq:pinn-bounded-H},
\begin{equation}\label{eq:pinn-Q-bound}
    \norm{\mathcal Q_k}_{L^2(Q_T)}
    \le
    \frac{M_2}{2}
    \norm{|Dh^k|^2}_{L^2(Q_T)} .
\end{equation}
\end{lem}

\begin{proof}
Using the definition of \(f\), we compute
\begin{align*}
    f[\widehat u^{k+1}]
    &=
    f[\widehat u^k+h^k] \\
    &=
    f[\widehat u^k]
    +(h^k)_t
    +D_pH(x,D\widehat u^k)\cdot Dh^k
    -\Delta h^k \\
    &\quad
    +H(x,D\widehat u^k+Dh^k)
    -H(x,D\widehat u^k)
    -D_pH(x,D\widehat u^k)\cdot Dh^k \\
    &=
    \rho_k+\mathcal Q_k .
\end{align*}
This proves \eqref{eq:pinn-residual-identity}.  Taylor's theorem in the
gradient variable gives
\[
    \mathcal Q_k
    =
    \int_0^1(1-s)
    D^2_{pp}H(x,D\widehat u^k+sDh^k)[Dh^k,Dh^k]\,ds .
\]
Therefore
\[
    |\mathcal Q_k|
    \le
    \frac{M_2}{2}|Dh^k|^2,
\]
and taking the \(L^2(Q_T)\)-norm gives \eqref{eq:pinn-Q-bound}.
\end{proof}

\begin{thm}[A posteriori PINN residual certificate]
\label{thm:pinn-aposteriori}
Let \(u\) be the classical periodic solution of \eqref{eq:pinn-hj}.  Let
\(\widehat u^k\) and \(h^k\) be sufficiently smooth periodic functions
satisfying
\[
    \widehat u^k(x,0)=g(x),
    \qquad
    h^k(x,0)=0
    \qquad\text{on } \T^n .
\]
Set
\[
    \widehat u^{k+1}:=\widehat u^k+h^k,
\]
and let \(\rho_k\) be the linearized Newton residual defined by
\eqref{eq:pinn-rho-natural}.  Assume \eqref{eq:pinn-bounded-H}.  Then
\begin{equation}\label{eq:pinn-aposteriori}
    \norm{\widehat u^{k+1}-u}_{\cE(Q_T)}
    \le
    C(T,M_1)
    \left(
        \norm{\rho_k}_{L^2(Q_T)}
        +
        \frac{M_2}{2}
        \norm{|Dh^k|^2}_{L^2(Q_T)}
    \right).
\end{equation}
Equivalently, defining the residual certificate indicator
\begin{equation}\label{eq:pinn-certificate}
    \mathfrak C_{k+1}
    :=
    \norm{\rho_k}_{L^2(Q_T)}
    +
    \frac{M_2}{2}
    \norm{|Dh^k|^2}_{L^2(Q_T)},
\end{equation}
we have
\begin{equation}\label{eq:pinn-certificate-short}
    \norm{\widehat u^{k+1}-u}_{\cE(Q_T)}
    \le
    C(T,M_1)\mathfrak C_{k+1}.
\end{equation}
\end{thm}

\begin{proof}
By Lemma \ref{lem:pinn-residual-identity},
\[
    f[\widehat u^{k+1}]
    =
    \rho_k+\mathcal Q_k .
\]
Using \eqref{eq:pinn-Q-bound}, we obtain
\[
    \norm{f[\widehat u^{k+1}]}_{L^2(Q_T)}
    \le
    \norm{\rho_k}_{L^2(Q_T)}
    +
    \frac{M_2}{2}
    \norm{|Dh^k|^2}_{L^2(Q_T)} .
\]
Also,
\[
    \widehat u^{k+1}(x,0)
    =
    \widehat u^k(x,0)+h^k(x,0)
    =
    g(x).
\]
Therefore Lemma \ref{lem:pinn-nonlinear-stability}, applied with
\(v=\widehat u^{k+1}\), gives
\[
    \norm{\widehat u^{k+1}-u}_{\cE(Q_T)}
    \le
    C(T,M_1)
    \norm{f[\widehat u^{k+1}]}_{L^2(Q_T)} .
\]

\end{proof}

\begin{rem}
The certificate in Theorem \ref{thm:pinn-aposteriori} is a deterministic
estimate in terms of continuous residuals.  
In computations, these residuals
are approximated by quadrature or sampling, so the sampling error is separate
from the PDE stability estimate proved here.  
Previous training residuals do
not appear explicitly in \eqref{eq:pinn-certificate-short}: they are already
encoded in the current iterate \(\widehat u^k\).  The updated iterate
\(\widehat u^{k+1}\) is verified through its current nonlinear residual,
expressed through the linearized residual \(\rho_k\) and the quadratic
Newton remainder.  The second term can also be written as
\[
    \norm{|Dh^k|^2}_{L^2(Q_T)}
    =
    \norm{Dh^k}_{L^4(Q_T)}^2 .
\]
\end{rem}

\subsection{Algorithmic form}

We summarize the Newton--PINN procedure and its residual certificate.
After \(K\ge1\) inexact Newton steps, Theorem~\ref{thm:pinn-aposteriori}
gives
\[
    \norm{\widehat u^K-u}_{\cE(Q_T)}
    \le C(T,M_1)\mathfrak C_K.
\]
Section~\ref{sec:numerics} examines this estimate numerically.
When an exact or manufactured solution \(u_{\rm ex}=u\) of
\eqref{eq:pinn-hj} is available, we compare \(\mathfrak C_K\) with
the true error
\[
    E_K:=\norm{\widehat u^K-u_{\rm ex}}_{\cE(Q_T)}
\]
through the calibration ratio \(E_K/\mathfrak C_K\), defined when
\(\mathfrak C_K>0\). By \eqref{eq:pinn-certificate-short}, this ratio
is bounded above by \(C(T,M_1)\).
A ratio that remains approximately constant and bounded away from zero
as training accuracy improves indicates that the indicator tracks
changes in the true error over the tested runs. Where no independent
solution is available, \(\mathfrak C_K\) remains a reference-free error
indicator; the error upper bound is \(C(T,M_1)\mathfrak C_K\).

\begin{algorithm}[t]
\caption{Newton--PINN iteration with residual certificate}
\label{alg:pinn-newton-certified}
\begin{algorithmic}[1]

\Require Hamiltonian \(H\) satisfying \eqref{eq:pinn-bounded-H},
initial data \(g\), final time \(T>0\), and an integer \(K\ge1\), the number of Newton steps.
\State Choose a smooth periodic \(\widehat u^0\) satisfying
\(\widehat u^0(x,0)=g(x)\).
\For{\(k=0,\dots,K-1\)}
    \State Compute
    \[
        \widehat R_k
        =
        \partial_t\widehat u^k
        +
        H(x,D\widehat u^k)
        -
        \Delta\widehat u^k .
    \]
    \State Freeze the Newton drift
    \[
        \widehat b^k
        =
        D_pH(x,D\widehat u^k).
    \]
\State Train a smooth periodic correction \(h^k\), with \(h^k(x,0)=0\),
to minimize a sampled approximation of
\(\norm{\rho_k}_{L^2(Q_T)}^2\), where
    \[
        \rho_k
        =
        (h^k)_t
        +
        \widehat b^k\cdot Dh^k
        -
        \Delta h^k
        +
        \widehat R_k .
    \]
    \State Update
    \[
        \widehat u^{k+1}
        =
        \widehat u^k+h^k .
    \]
    \State Record
    \[
        \mathfrak C_{k+1}
        =
        \norm{\rho_k}_{L^2(Q_T)}
        +
        \frac{M_2}{2}
        \norm{|Dh^k|^2}_{L^2(Q_T)} .
    \]
\EndFor
\Ensure Approximation \(\widehat u^K\) and certificate indicator
\(\mathfrak C_K\).
\end{algorithmic}
\end{algorithm}

% Include a figure if the PDF exists; otherwise show a placeholder box.
% (Sweep figures are generated by code_final/make_paper_figures.py.)
\providecommand{\maybegraphics}[2]{%
  \IfFileExists{#2}{\includegraphics[#1]{#2}}{\fbox{\ttfamily pending: \detokenize{#2}}}}

\section{Numerical experiments}
\label{sec:numerics}

We examine the accuracy and convergence behavior of inexact Newton--PINN
iterations, the effect of frozen linearization on training cost, and
the usefulness of the residual indicator for monitoring and stopping.
A one-dimensional benchmark with an exact solution
is followed by nonconvex high-dimensional problems and comparisons that separate network capacity, staged training, and linearization.
We then study per-step training budgets and fixed-budget allocation across Newton steps, evaluate the indicator on wide-interaction and unforced high-dimensional benchmarks---including regimes of failure and stagnation---and test indicator-driven stopping.

\subsection{Common method, error measures, and experimental protocol}
\label{subsec:protocol}

\subsubsection{Iteration and residual indicator.}
Throughout the experiments, $\T^n=(\mathbb R/(2\pi\mathbb Z))^n$
and $Q_T=\T^n\times(0,T)$.
For a possibly zero manufactured source $S$, write
\[
    \mathcal R[v]:=v_t+H(Dv)-\Delta v-S.
\]
The proof of Theorem~\ref{thm:pinn-aposteriori} applies with the same
prescribed source $S$ in both equations, since $S$ cancels in the
residual difference.
The approximation has the form
\[
    \widehat u^K(x,t)=g(x)+\sum_{j=0}^{K-1}h^j(x,t),
    \qquad
    h^j(x,t)=t\,\mathcal N^j(\sin x,\cos x,t/T),
\]
where $\mathcal N^j$ is a fully connected $\tanh$ network and the spatial
features are coordinatewise. 
For periodic $g$, this ansatz satisfies the periodic boundary
and initial conditions exactly.
 At step $k$, a fresh network is trained with
$D_pH(D\widehat u^k)$ and $\mathcal R[\widehat u^k]$ frozen, and
$\widehat u^{k+1}=\widehat u^k+h^k$, $k=0,\ldots,K-1$.
The linearized residual is
\[
    \rho_k=(h^k)_t+D_pH(D\widehat u^k)\cdot Dh^k
           -\Delta h^k+\mathcal R[\widehat u^k].
\]
The new network is trained on the sampled squared linearized residual,
with the normalization specified below.
Identity~\eqref{eq:pinn-residual-identity} decomposes the updated nonlinear
residual into this term and a quadratic Taylor remainder. The resulting
continuous indicator is
\begin{equation}\label{eq:numerical-certificate}
    \mathfrak C_K
    :=\norm{\rho_{K-1}}+
      \frac{M_2}{2}\norm{|Dh^{K-1}|^2}, \qquad K\ge1.
\end{equation}
The norms are the $L^2(Q_T)$ norms of
Section~\ref{sec:pinn-certificate} in Experiments~1 and~3 and their
root-mean-square (RMS) normalized counterparts in the high-dimensional
experiments.

Quadrature or Monte Carlo sampling gives an estimate
$\widehat{\mathfrak C}_K$, 
denoted by $\mathfrak C_K$ when numerical values are reported.
Theorem~\ref{thm:pinn-aposteriori} bounds the energy error by
$C(T,M_1)\mathfrak C_K$ under its hypotheses. We evaluate neither this
constant nor a confidence bound for quadrature error, so the numerical
indicator is used empirically: a tolerance is an indicator threshold,
and calibration refers to $E_K/\widehat{\mathfrak C}_K$ when the exact
solution is available.

\subsubsection{Architecture and training.}
All runs use $T=1$ and networks with five hidden layers, of width $96$
in Experiments~1 and~3 and $128$ in Experiment~2. Adam uses learning rate $10^{-3}$ and fresh uniform collocation batches; Experiment~1 adds an
L-BFGS refinement after each Adam phase. In the high-dimensional Newton--PINN experiments, each step's training loss is divided by a fixed positive estimate of the squared RMS of $\mathcal R[\widehat u^k]$. The denominator is held fixed during each inner solve, so this rescaling preserves the minimizer; reported residuals are not rescaled. Experiments~1 and~3 use the unscaled sampled squared linearized residual.

 Training batches contain $4096$ points for $n\le20$, $2048$
for $n=50$, and $1024$ for $n=100$. These reductions keep the
automatic-differentiation Laplacian within GPU memory, but change the
work per epoch; reported times therefore describe this memory-adapted
protocol rather than fixed-work scaling. For Adam, one epoch denotes one optimizer update using a fresh collocation batch.

\subsubsection{Error measurement.}
We denote the energy error by $E_K$. Experiments~1 and~3 use
\eqref{eq:pinn-energy-norm}; high-dimensional experiments use
\[
\begin{split}
    E_K^2:={}&\max_{t\in[0,T]}\frac{1}{(2\pi)^n}
       \int_{\T^n}|\widehat u^K-u_{\rm ex}|^2\,dx\\
    &+\frac{1}{T(2\pi)^n}\int_0^T\!\!\int_{\T^n}
       |D\widehat u^K-Du_{\rm ex}|^2\,dx\,dt,
\end{split}
\]
with residual norms normalized by
\[
    \norm{w}_{\rm RMS}^2
    :=\frac{1}{|Q_T|}\int_{Q_T}|w|^2\,dx\,dt.
\]
This removes volume factors from dimension comparisons. In Experiment~1, the energy error is evaluated on a uniform grid
with $512$ spatial points and $241$ time levels, using trapezoidal integration in time and periodic central differences for the gradient term. Experiment~2 uses $8\times10^4$, $4\times10^4$, and $2\times10^4$
uniform space--time Monte Carlo samples for each residual and
indicator evaluation at $n\le20$, $n=50$, and $n=100$, respectively.
The energy error $E_K$ uses $31$ uniform time slices, with
$8192$, $4096$, and $2048$ fresh spatial samples per slice,
respectively. Error, residual, and indicator evaluations use separate samples.

To assess evaluation error, one trained model from each regime
(well-resolved, $n=10$; training floor, $n=50$; partial resolution,
$n=100$; failed wide family, $n=20$) is re-evaluated under ten independent evaluation
seeds at each $N\in\{2\times10^4,8\times10^4,3.2\times10^5\}$.
At the tested sample sizes, relative standard deviations are at most
$1.4\%$ for error and $0.8\%$ for indicator and residual; 
the mean estimates differ by at most $1\%$ from those obtained
with $N=3.2\times10^5$.
For these four models, evaluation variability is small relative to the
large error reductions.

All runs use a single NVIDIA GeForce RTX~5070 GPU (12\,GB) and
PyTorch~2.11. Wall-clock training times and peak CUDA memory are 
reported with the aggregation specified in each caption, excluding evaluation and figure generation.

\subsection{Baselines and relation to neural HJB solvers}
\label{subsec:numerics-related}

Neural PDE solvers include direct physics-informed residual minimization
\cite{raissi2019physics}. Related methods for HJB equations include
adaptive data generation \cite{NakamuraZimmererGongKang2021}, actor--critic methods for static HJB equations \cite{ZhouHanLu2021}, trajectory-based training for stochastic control
\cite{LiVermaRuthotto2024}, and randomized-smoothing PINNs designed to reduce
the cost of derivative evaluation in high dimensions
\cite{HuYangWangKarniadakisKawaguchi2025}. Our comparisons examine
Newton-based training for viscous Hamilton--Jacobi equations with possibly
nonconvex Hamiltonians; they do not rank these other neural solvers.

The direct PINN baseline (C1) uses the ansatz
\[
    v_\theta(x,t)=g(x)+t\,\mathcal N_\theta(\sin x,\cos x,t/T),
\]
where $\mathcal N_\theta$ is a single network with the same width and
depth as a Newton correction network. It is trained on a
sampled approximation of $\|\mathcal R[v_\theta]\|_{\rm RMS}^2$.
Since Newton accumulates correction networks, we report both training
budgets and parameter counts. Section~\ref{subsec:exp-baselines} adds
composite and sequential nonlinear baselines to separate capacity,
staging, and linearization. 
The manufactured high-dimensional tests assess solver accuracy and
residual monitoring rather than performance on a specific
stochastic-control application.

\subsection{Experiment 1: exact Cole--Hopf benchmark and initial contraction}
\label{subsec:exp-cole-hopf}

The first experiment checks accuracy against an exact solution and
examines the transition from rapid improvement to a training plateau.
Consider
\begin{equation}\label{eq:exp-cole-hopf-pde}
\begin{cases}
    u_t+\dfrac12|u_x|^2=u_{xx}
       \qquad &\text{in } \T\times(0,1),\\[1mm]
    u(x,0)=g(x) \qquad &\text{on } \T,
\end{cases}
\end{equation}
with $H(p)=p^2/2$, so $M_2=1$.
Since $D_pH(p)=p$ is not globally bounded, 
the stability argument uses a bound on the mean-value drift
$(u_x+(\widehat u^K)_x)/2$.
 For
$\phi(x,t)=1+ae^{-t}\cos x$, $0<a<1$, the Cole--Hopf solution is
\[
    u_{\rm ex}(x,t)=-2\log\phi(x,t),
    \qquad g(x)=-2\log(1+a\cos x).
\]
We take $a=1/2$ and use double precision, $K=5$ corrections, $4000$ Adam
epochs followed by $400$ strong-Wolfe L-BFGS iterations per correction,
and batches of $4096$. A semi-implicit Fourier finite-difference method
(FDM), with $512$ spatial grid points and $4000$ time steps, provides a
classical reference.

\begin{figure}[htbp]
    \centering
    \includegraphics[width=\textwidth]{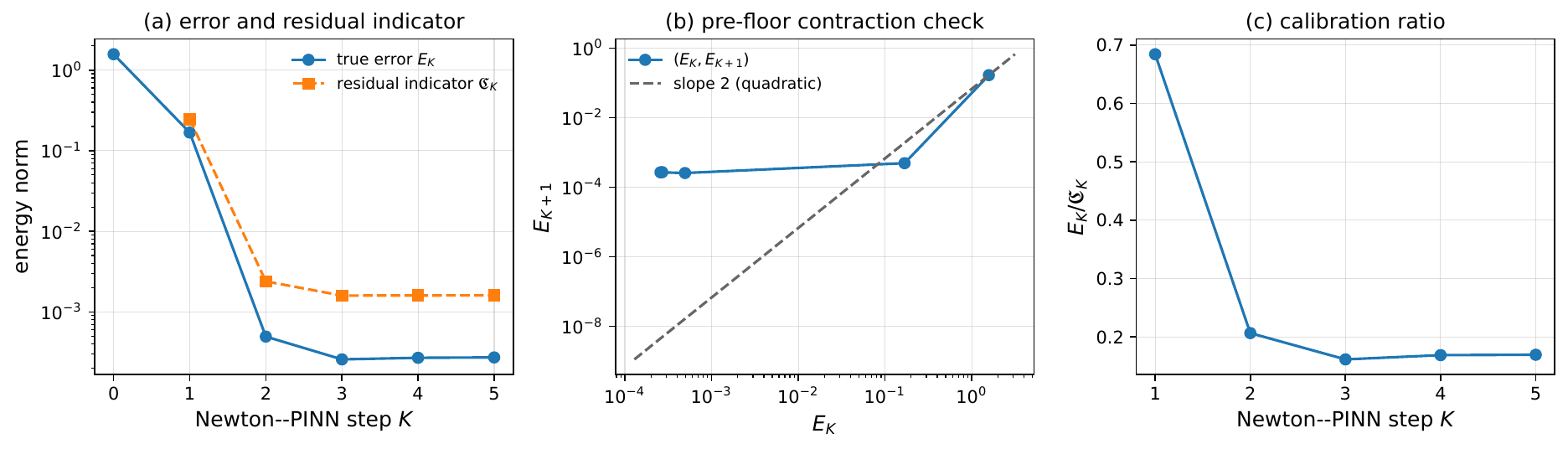}
    \caption{Experiment 1 (Cole--Hopf, $n=1$).
    (a) Energy error $E_K$ and empirical indicator $\mathfrak C_K$.
    (b) Pairs $(E_K,E_{K+1})$ and the slope-$2$ reference through the
    first pair; the transition $E_1\to E_2$ is consistent with the local
    quadratic mechanism. (c) Calibration ratio $E_K/\mathfrak C_K$.}
    \label{fig:exp1-main}
\end{figure}

\subsubsection{Accuracy and plateau.}
Figure~\ref{fig:exp1-main} gives
\[
\begin{aligned}
 E_0&=1.57,& E_1&=1.68\times10^{-1},& E_2&=4.94\times10^{-4},\\
 E_3&=2.58\times10^{-4},& E_4&=2.70\times10^{-4},& E_5&=2.72\times10^{-4}.
\end{aligned}
\]
The pair $(E_1,E_2)$ lies below the slope-$2$ line through $(E_0,E_1)$.
This is a consistency check on the local quadratic mechanism: only one
informative transition occurs before the plateau, and the energy norm
differs from the Schauder norm of Theorem~\ref{thm:main2}. Three
corrections reach accuracy comparable to the chosen FDM reference
($1.57\times10^{-4}$ on the $512\times4000$ grid); further corrections
do not reduce the error at the prescribed per-step budget.

\subsubsection{Indicator behavior.}
The indicator decreases with the error:
$\mathfrak C_1=2.45\times10^{-1}$,
$\mathfrak C_2=2.39\times10^{-3}$, and
$\mathfrak C_3=1.60\times10^{-3}$.
For $K\ge3$, $E_K/\mathfrak C_K$ stabilizes at $0.16$--$0.17$, so the
unscaled indicator is about six times the error at the plateau.
The quadratic term $\frac{M_2}{2}\norm{|Dh^0|^2}$ dominates at $K=1$,
and $\norm{\rho_{K-1}}$ dominates from $K=2$ onward, illustrating both
terms of \eqref{eq:numerical-certificate}.

\subsection{Experiment 2: nonseparable nonconvex problems in dimensions 10--100}
\label{subsec:exp-highdim}

This experiment tests accuracy and residual monitoring for nonseparable,
nonconvex problems in high dimensions. 
To evaluate accuracy without a high-dimensional grid reference,
we use manufactured solutions on
$\T^n\times[0,1]$:
\begin{equation}\label{eq:exp-hd-pde}
\begin{cases}
    u_t+H(Du)-\Delta u=S(x,t),\\
    u(x,0)=u_{\rm ex}(x,0),
\end{cases}
\end{equation}
where $S:=(u_{\rm ex})_t+H(Du_{\rm ex})-\Delta u_{\rm ex}$ is computed
analytically. 
Pointwise solution and gradient errors can therefore be evaluated
against analytic references in every tested dimension.
Before training, the closed-form derivatives of every instance
are checked against double-precision automatic differentiation, with
agreement to $10^{-13}$.

The nonconvex Hamiltonian couples all gradient components:
\begin{equation}\label{eq:exp-hd-H}
    H(p)
    =
    \frac1 n\sum_{i=1}^n
    \Bigl(\sin p_i+\frac14\cos(2p_i)\Bigr)
    +\beta\sin\Bigl(\frac{1}{\sqrt n}\sum_{i=1}^n p_i\Bigr)
    +\gamma\cos(w\cdot p),
\end{equation}
with $\beta=0.35$, $\gamma=0.25$, and $w=(e_1-e_3+e_5)/\sqrt3$.
In \eqref{eq:numerical-certificate}, we replace $M_2$ by the explicit
upper bound $\overline M_2$:
\[
    \overline M_2:=\frac2n+|\beta|+|\gamma|\ge M_2.
\]
The quadratic indicator coefficient $\overline M_2/2$
tends to $0.3$ as $n\to\infty$.

Set
\[
    P(x)=\frac{1}{\sqrt n}\sum_{i=1}^n(-1)^{i-1}\cos(x_i+x_{i+1}),
    \qquad
    Q(x)=\frac{1}{\sqrt n}\sum_{i=1}^n(-1)^{i-1}\sin(x_i+x_{i+2}),
\]
with cyclic indices, so that $\Delta P=-2P$ and $\Delta Q=-2Q$.
The two solution families are
\begin{align}
    u^{\rm chain}_{\rm ex}(x,t)
    &=
    \sin\Bigl(\sum_{i=1}^{5}x_i\Bigr)
    +t\,P(x)+e^{-t}Q(x),
    \label{eq:chain-family}\\
    u^{\rm ring}_{\rm ex}(x,t)
    &=
    \sum_{j=1}^{n}f_j\bigl(x_j+2x_{j+1}-x_{j+2}\bigr)
    +t\,P(x)+e^{-t}Q(x),
    \label{eq:ring-family}
\end{align}
where $f_j=\sin$ for odd $j$ and $f_j=\cos$ for even $j$. Both families
are $2\pi$-periodic and nonseparable. Their common time-dependent part
involves every coordinate and has RMS scale $O(1)$ uniformly in $n$.
The time-independent background involves only the first five coordinates
in the chain family and all coordinates through overlapping triples in
the ring family. In both cases, the frozen drift $D_pH(D\widehat u^k)$ generally has nonzero components in every coordinate direction. The $O(1)$ signal scale and nonvanishing global Hessian bound do not, by themselves, establish dimension-independent nonlinear difficulty.

\subsubsection{Protocol and comparisons.}
The main sweep uses single precision and $K=4$ Newton corrections, each
trained for $4000$ Adam epochs. The direct PINN baseline (C1) is trained for $16000$ Adam epochs. Thus the full $K=4$ comparison
matches epochs, but the direct baseline has one quarter of the Newton
composite's parameters. Section~\ref{subsec:exp-baselines} separates
capacity, staging, and linearization; its matched comparisons are
essential to interpreting the baseline results below.

The coordinatewise automatic-differentiation Laplacian dominates the cost
per epoch. Both families use three random seeds at $n\in\{10,20\}$;
at $n\in\{50,100\}$, the ring family uses five seeds and the chain family
one. Supplementary ring runs vary the number of corrections and the
per-step budget separately; their protocols and results appear in
Section~\ref{subsec:budget-allocation}.

\begin{figure}[htbp]
    \centering
    \maybegraphics{width=\textwidth}{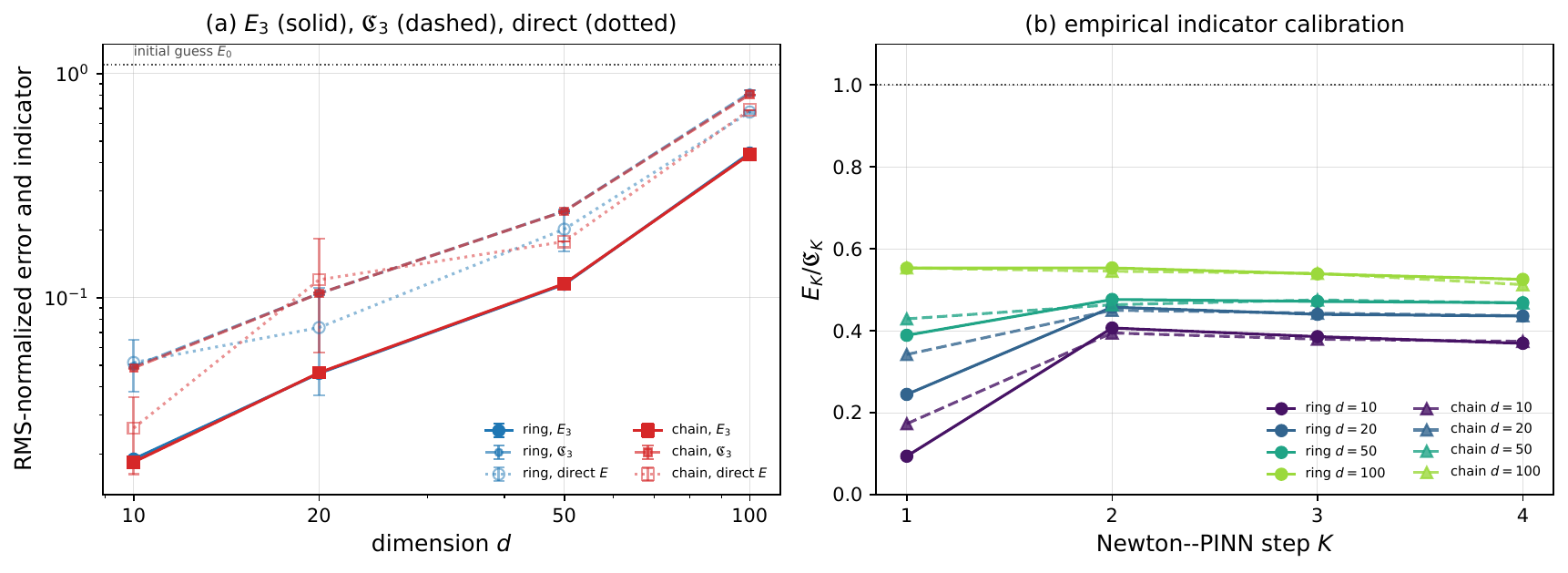}
    \caption{Experiment 2 (dimensions $10$--$100$). (a) Seed-averaged
    energy error $E_3$ (solid) and empirical residual indicator
    $\mathfrak C_3$ (dashed) versus dimension for both families, with the
    direct PINN after the full $K=4$ epoch budget (dotted) and initial
    error $E_0$ (horizontal). Error bars denote sample standard deviations where at least two seeds
are available.
    (b) Empirical calibration ratio $E_K/\mathfrak C_K$ versus $K$ for
    both families.}
    \label{fig:sweep-main}
\end{figure}

\begin{table}[htbp]
    \centering
    \small
    \setlength{\tabcolsep}{5pt}
    \resizebox{\textwidth}{!}{% Auto-generated by make_paper_figures.py -- do not edit by hand.
\begin{tabular}{lrcccccccc}
\toprule
family & $n$ & $E_{3}$ & $\mathfrak{C}_{3}$ & $E_{3}/\mathfrak{C}_{3}$ & $\|\mathcal R[\widehat u^{3}]\|$ & $E^{\rm direct}$ & $\|\mathcal R[v^{\rm direct}]\|$ & $t_{K=3}$ (s) & mem (GB) \\
\midrule
chain & 10 & $1.85(2){\times}10^{-2}$ & $4.87(5){\times}10^{-2}$ & $0.380(5)$ & $4.88(6){\times}10^{-2}$ & $2.62(99){\times}10^{-2}$ & $4.46(67){\times}10^{-2}$ & 326 & 1.33 \\
chain & 20 & $4.63(7){\times}10^{-2}$ & $1.05(1){\times}10^{-1}$ & $0.443(2)$ & $1.04(1){\times}10^{-1}$ & $1.20(63){\times}10^{-1}$ & $1.55(52){\times}10^{-1}$ & 582 & 2.63 \\
chain & 50 & $1.16{\times}10^{-1}$ & $2.43{\times}10^{-1}$ & $0.48$ & $2.44{\times}10^{-1}$ & $1.78{\times}10^{-1}$ & $3.02{\times}10^{-1}$ & 1220 & 3.52 \\
chain & 100 & $4.36{\times}10^{-1}$ & $8.08{\times}10^{-1}$ & $0.54$ & $7.88{\times}10^{-1}$ & $6.89{\times}10^{-1}$ & $1.15{\times}10^{0}$ & 2398 & 4.02 \\
\midrule
ring & 10 & $1.91(7){\times}10^{-2}$ & $4.94(10){\times}10^{-2}$ & $0.39(1)$ & $4.94(12){\times}10^{-2}$ & $5.1(13){\times}10^{-2}$ & $7.2(20){\times}10^{-2}$ & 329 & 1.33 \\
ring & 20 & $4.60(10){\times}10^{-2}$ & $1.04(2){\times}10^{-1}$ & $0.440(4)$ & $1.04(2){\times}10^{-1}$ & $7.4(37){\times}10^{-2}$ & $1.12(27){\times}10^{-1}$ & 615 & 2.64 \\
ring & 50 & $1.15(2){\times}10^{-1}$ & $2.43(4){\times}10^{-1}$ & $0.472(3)$ & $2.44(3){\times}10^{-1}$ & $2.02(41){\times}10^{-1}$ & $3.19(30){\times}10^{-1}$ & 1345 & 3.56 \\
ring & 100 & $4.42(16){\times}10^{-1}$ & $8.20(23){\times}10^{-1}$ & $0.54(1)$ & $7.97(18){\times}10^{-1}$ & $6.74(9){\times}10^{-1}$ & $1.13(1){\times}10^{0}$ & 2662 & 4.09 \\
\bottomrule
\end{tabular}}
    \caption{Experiment 2: RMS-normalized energy error, empirical
    indicator, calibration ratio, Newton nonlinear residual
    $\norm{\mathcal R[\widehat u^3]}_{\rm RMS}$, direct-PINN error after
    $16000$ epochs and its nonlinear residual
    $\norm{\mathcal R[v^{\rm direct}]}_{\rm RMS}$, cumulative Newton
    training time to $K=3$, and peak CUDA memory over the full $K=4$ run.
    The residual columns compare indicators governed by the same
    stability estimate. Parenthesized digits give the sample standard
    deviation over seeds in the last quoted place; single-seed entries
    have none.}
    \label{tab:sweep}
\end{table}

\begin{figure}[htbp]
    \centering
    \maybegraphics{width=0.62\textwidth}{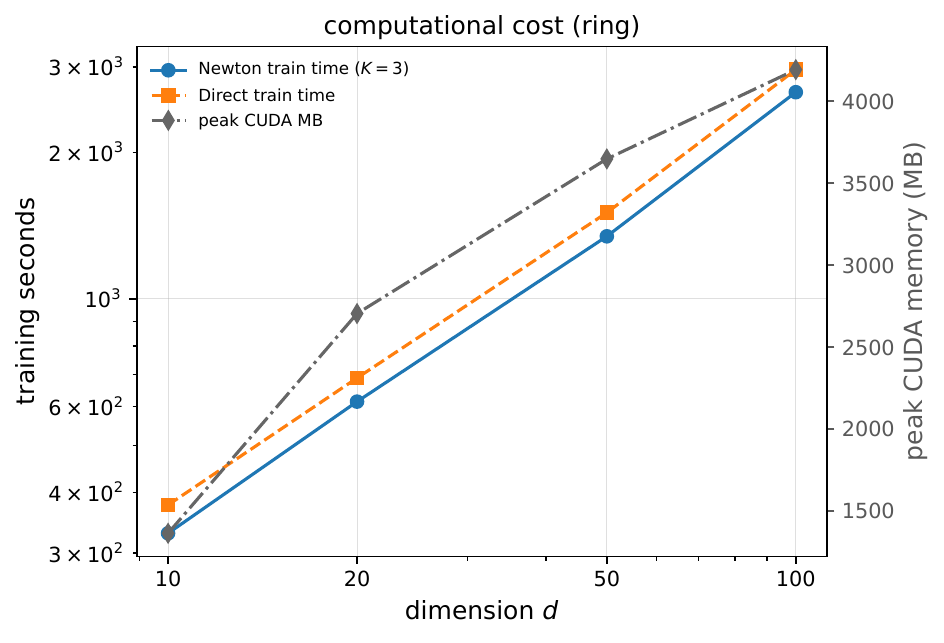}
    \caption{Experiment 2, computational cost (ring family): cumulative
    Newton training time to $K=3$, direct-baseline training time for $16000$ epochs, and peak CUDA memory over the full $K=4$ run, versus dimension under the memory-adapted batch protocol of Section~\ref{subsec:protocol}.}
    \label{fig:sweep-cost}
\end{figure}

Figure~\ref{fig:sweep-main} and Table~\ref{tab:sweep} collect the main
results; Figures~\ref{fig:sweep-budget} and~\ref{fig:sweep-cost} isolate
budget dependence and cost. Two-dimensional slices are supplementary
diagnostics: fixing the other $n-2$ coordinates changes the signal scale
with $n$, so these plots alone do not support cross-dimensional
accuracy comparisons.

\subsubsection{Accuracy relative to the learned increment.}
Because the ansatz imposes the initial condition exactly, the corrections
are trained to approximate $u_{\rm ex}(x,t)-u_{\rm ex}(x,0)$. This increment is the same for both
families, has RMS scale $O(1)$ uniformly in $n$, and has energy norm
$E_0\approx1.10$. The iteration reduces this initial error by $98\%$ at
$n=10$, $96\%$ at $n=20$, and $89\%$ at $n=50$ ($92\%$ with the doubled
per-step budget). At $n=100$, the reduction is $59\%$, rising to $78\%$
at $K=8$ with the doubled per-step budget. These $n=100$ iterates only
partially resolve the solution; the indicator reports this without using
$u_{\rm ex}$.

\subsubsection{Dimension dependence.}
The families' energy errors $E_3$ differ by at most $4\%$ at the tested
dimensions (Table~\ref{tab:sweep}). This comparison concerns different
known spatial backgrounds with the same increment to be learned.
At the standard epoch count and under the dimension-dependent batch
protocol, the ring family's log--log slope of $E_3$ versus $n$
is $1.12$ from $n=10$ to $n=50$ and $1.95$ from $n=50$ to $n=100$. Doubling the per-step budget at both endpoints lowers the latter slope to $1.83$. These slopes are budget-sensitive and do not establish an asymptotic dimension law; we therefore report no single fitted exponent. The budget study in Section~\ref{subsec:budget-allocation} examines this sensitivity further.

The energy norm includes the gradient and the $\max_t$ spatial norm. In the same RMS normalization, the space-time $L^2(Q_T)$ error alone is
$0.38$--$0.44$ times the energy error for the seed-averaged $K=3$ Newton results in the main sweep (Table~\ref{tab:sweep}). At the doubled per-step budget, the two-seed ring averages of $\norm{\widehat u^K-u_{\rm ex}}_{\rm RMS}$ are $3.7\times10^{-2}$ at $n=50$ ($K=3$) and
$9.9\times10^{-2}$ at $n=100$ ($K=8$), compared with $9.9\times10^{-2}$ and $3.0\times10^{-1}$ for the direct baseline after $16000$ and $64000$ epochs, respectively. These are absolute errors against an $O(1)$ time-dependent signal, not relative errors.

\subsubsection{Indicator calibration.}
Across the main-sweep families, dimensions, and seeds, $E_K/\mathfrak C_K$ is approximately $0.37$--$0.50$ for $n\le50$ and $0.51$--$0.56$ for $n=100$ when $2\le K\le4$ (Figure~\ref{fig:sweep-main}(b)). This band is narrower than a factor $1.6$ while the error spans a factor $33$: the unscaled indicator overestimates the measured error by factors $1.8$--$2.7$. The doubled-budget $n=50$ runs remain in the same band ($0.46$--$0.48$). After the first correction, the quadratic term
$\frac{\overline M_2}{2}\norm{|Dh^0|^2}$ dominates at $n=10$, as in Experiment~1. The first-correction calibration ratios span $0.08$--$0.45$ for $n\le50$. At $n=100$, the first correction removes only about $30\%$ of the error, the quadratic term does not dominate, and the ratio is $0.55$--$0.56$.

\subsubsection{Residuals and the direct baseline.}
The stability estimate behind Theorem~\ref{thm:pinn-aposteriori} applies
to any sufficiently smooth periodic approximation satisfying the initial condition exactly. A direct
PINN therefore has the analogous indicator $\norm{\mathcal R[v]}$;
residual-based error control is not specific to Newton. Accordingly,
Table~\ref{tab:sweep} reports both methods' nonlinear residuals. The reported ratio $\norm{\mathcal R[\widehat u^4]}/\mathfrak C_4$ is $0.975$--$1.005$ over all $33$ runs of the sweep and wide family: the indicator from the final linear solve and the directly evaluated nonlinear residual are evaluated on independent Monte Carlo samples and are numerically close. The timing logs do not separate their evaluation costs, so this shows agreement, not a speed advantage.

At $K=4$ and $n\ge50$, the direct baseline's seed-averaged residual is $1.3$--$1.6$ times Newton's (five ring seeds and one chain seed at each dimension), compared with $1.2$--$1.5$ at the $K=3$ iterate in Table~\ref{tab:sweep}.

For $n\le20$ at $K=4$, the residuals are comparable and their
ordering depends on the seed (ratios $0.85$--$1.97$). Newton then has a smaller and less variable true error, but not consistently a smaller indicator.

At the $K=3$ comparison, Newton has used $12000$ epochs against the direct
baseline's $16000$. In the ring family, the baseline's errors are larger
by factors $2.70$, $1.60$, $1.76$, and $1.52$ at $n=10,20,50,100$,
respectively (three seeds at $n\le20$, five at $n\ge50$).
The largest factor occurs at the smallest dimension.
The direct baseline's seed spread is $26$--$53\%$ of its mean
across both families at $n\le20$ and $20\%$ for the ring family
at $n=50$, compared with at most $4\%$ for Newton ($1.4\%$ and $3.5\%$ for the ring family at $n=50,100$). This comparison gives the direct
baseline more epochs but controls neither wall time nor parameter count;
Section~\ref{subsec:exp-baselines} separates these effects at $K=4$.

\subsubsection{Computational cost.}
Under the memory-adapted protocol, training to $K=3$ takes $5.5$--$45$\,min
and peak memory rises from $1.3$ to $4.1$\,GB as $n$ grows from $10$ to
$100$. On the only segment with a fixed training batch, $n=10\to20$,
the ring training time rises from $329$ to $615$\,s. From $n=50$ to
$n=100$, time nearly doubles while the batch halves, making the time
per sampled collocation point about four times larger. These results
document single-GPU feasibility through $n=100$, not sublinear scaling
with dimension.

\subsection{Separating capacity, staging, and linearization}
\label{subsec:exp-baselines}

 Two additional baselines on the ring family
at $n\in\{50,100\}$, with two seeds each, separate parameter count,
staging, and linearization. The composite direct baseline (C2) trains
the identical architecture $g+t\sum_{j=1}^{4}\mathcal N_j$ jointly on
the nonlinear residual for $16000$ epochs: equal capacity, without staging
or linearization. The sequential nonlinear baseline (C3) adds one network
at a time on the Newton schedule ($4\times4000$ epochs, the parameters of previous networks
held fixed), training each stage on the nonlinear residual of the whole
composite. It thus matches capacity growth and staging without Newton linearization.

\begin{table}[htbp]
    \centering
    \small
    \setlength{\tabcolsep}{5pt}
    % Auto-generated by make_phase123_tables.py -- do not edit by hand.
\begin{tabular}{llcccc}
\toprule
$n$ & method & params & $E$ & $\|\mathcal R\|$ & wall (s) \\
\midrule
50 & Newton--PINN ($K{=}4$) & $4N_{\rm par}$ & $1.13(1){\times}10^{-1}$ & $2.38{\times}10^{-1}$ & 2014 \\
50 & Newton--PINN ($K{=}4$, $8000$\,ep) & $4N_{\rm par}$ & $8.93(1){\times}10^{-2}$ & $1.90{\times}10^{-1}$ & 3971 \\
50 & Seq.\ nonlinear (C3, $K{=}4$) & $4N_{\rm par}$ & $1.12(1){\times}10^{-1}$ & $2.38{\times}10^{-1}$ & 7754 \\
50 & Composite direct (C2) & $4N_{\rm par}$ & $1.52(3){\times}10^{-1}$ & $2.15{\times}10^{-1}$ & 8942 \\
50 & Direct, single net (C1) & $N_{\rm par}$ & $1.93(14){\times}10^{-1}$ & $3.13{\times}10^{-1}$ & 1523 \\
\midrule
100 & Newton--PINN ($K{=}4$) & $4N_{\rm par}$ & $3.88(17){\times}10^{-1}$ & $7.20{\times}10^{-1}$ & 3991 \\
100 & Newton--PINN ($K{=}4$, $8000$\,ep) & $4N_{\rm par}$ & $2.77(9){\times}10^{-1}$ & $5.47{\times}10^{-1}$ & 9054 \\
100 & Seq.\ nonlinear (C3, $K{=}4$) & $4N_{\rm par}$ & $3.86(9){\times}10^{-1}$ & $7.32{\times}10^{-1}$ & 5681 \\
100 & Composite direct (C2) & $4N_{\rm par}$ & $3.35(20){\times}10^{-1}$ & $5.58{\times}10^{-1}$ & 8212 \\
100 & Direct, single net (C1) & $N_{\rm par}$ & $6.71(14){\times}10^{-1}$ & $1.13{\times}10^{0}$ & 2992 \\
\bottomrule
\end{tabular}
    \caption{Ring family, $n\in\{50,100\}$: energy error, nonlinear
    residual, and wall-clock training time for Newton and the
    capacity-matched baselines. Equal epoch counts need not give equal
    wall time because the differentiated computation graphs differ.
    $N_{\rm par}$ denotes one network's parameter count. Energy errors, residuals,
    and times are means over two seeds; parenthesized digits are as in
    Table~\ref{tab:sweep}.}
    \label{tab:baselines}
\end{table}

\subsubsection{Staging and linearization.}
C3 and standard Newton attain similar mean energy errors at both dimensions:
$1.12\times10^{-1}$ versus $1.13\times10^{-1}$ at $n=50$, and
$3.86\times10^{-1}$ versus $3.88\times10^{-1}$ at $n=100$
(Table~\ref{tab:baselines}). Thus the accuracy observed here is attainable
with staged representation growth without Newton linearization.
On the same $4\times4000$ epoch schedule, however, Newton takes
$2014$ versus $7754$\,s at $n=50$ and $3991$ versus $5681$\,s at $n=100$,
reducing training time by factors $3.9$ and $1.4$.

\subsubsection{Capacity and training time.}
At equal epoch and parameter counts, C2 reaches $1.52\times10^{-1}$
against Newton's $1.13\times10^{-1}$ at $n=50$; the ordering reverses
at $n=100$, with $3.35\times10^{-1}$ against $3.88\times10^{-1}$.
Epochs are not compute-matched: C2 differentiates through all four
networks at every epoch. At $n=50$, Newton with $8000$ epochs per step
achieves $8.93\times10^{-2}$ in $3971$\,s, 
giving a smaller error in less time than C2
 ($1.52\times10^{-1}$ in $8942$\,s).
At $n=100$, C2's error and training time lie between those of the
two Newton runs. Relative to C2, standard Newton uses
$51\%$ less time with $16\%$ larger error, while doubled-budget Newton
uses $10\%$ more time with $17\%$ smaller error. Neither $n=100$
comparison is exactly time-matched.

\subsection{Allocating the training budget}
\label{subsec:budget-allocation}

Two supplementary sets of ring-family runs, each using two seeds,
distinguish the benefit of more corrections from that of more training
per correction. The first extends the standard $4000$-epoch schedule
to $K=8$ at $n\in\{10,20,50,100\}$. The second uses $8000$ epochs per
step, to $K=8$ at $n=100$ and to $K=4$ at $n=50$.

\begin{figure}[htbp]
    \centering
    \maybegraphics{width=0.62\textwidth}{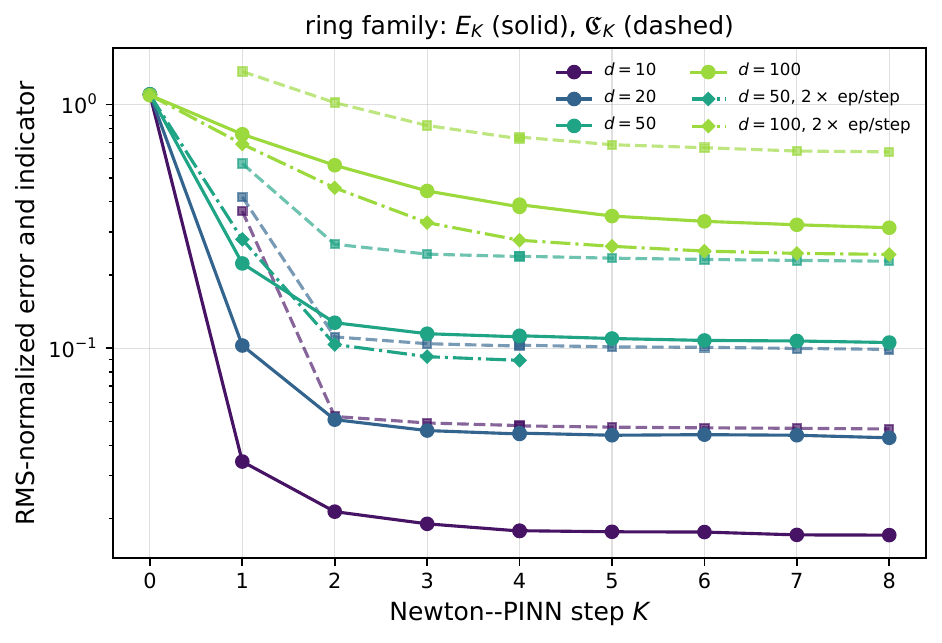}
    \caption{Experiment 2, saturation and budget dependence (ring family).
    Seed-averaged $E_K$ (solid) and $\mathfrak C_K$ (dashed) versus $K$.
    Supplementary runs 
    continue the standard-budget curves to $K=8$
     with two seeds and
    the same per-step budget. Dash-dotted curves use $8000$ epochs per
    step, to $K=8$ at $n=100$ and to $K=4$ at $n=50$.}
    \label{fig:sweep-budget}
\end{figure}

\subsubsection{Saturation and per-step training at $n\le50$.}
In every run, the error decreases monotonically with $K$ until it reaches
the training floor. The first correction reduces the $O(1)$ initial
error ($E_0\approx1.10$) by up to two orders of magnitude, and the
$n\le50$ runs saturate at $K=3$. Extending them to $K=8$ at the same
per-step budget changes the error by only $9\%$, $8\%$, and $9\%$
at $n=10,20,50$, respectively. The sequence is no longer monotone within
that range (Figure~\ref{fig:sweep-budget}), as in Experiment~1.

The plateau responds more to per-step training than to further corrections
in the tested $n=50$ runs. Doubling the epochs per step lowers $E_3$
from $1.16\times10^{-1}$ to $9.2\times10^{-2}$ and gives
$E_4=8.9\times10^{-2}$, reductions of $20\%$ and $23\%$ against the
$9\%$ obtained from five more corrections at the original budget.
The doubled-budget run stops at $K=4$: its per-step improvement has
fallen to $3\%$, compared with $15\%$ at $K=4$ for $n=100$.

\subsubsection{Continuing improvement at $n=100$.}
At $n=100$, the standard run has not saturated: the error falls by $14\%$
from $K=3$ to $K=4$, so $E_3$ is a budget-dependent snapshot rather than
a floor. Continuing at the same per-step budget gives
$E_8=3.12\times10^{-1}$, $31\%$ below $E_3$, while successive reductions
decline from $14\%$ to $3\%$. Doubling the per-step budget gives
$E_8=2.43\times10^{-1}$, $46\%$ below the standard-budget $E_3$
(Figure~\ref{fig:sweep-budget}). A pooled fit to the sixteen points
from both schedules over $B\in[4\times10^3,6.4\times10^4]$ total epochs
gives
\begin{equation}\label{eq:budget-law}
    E\;\propto\;B^{-0.46},
    \qquad
    R^2=0.96.
\end{equation}
The fit \eqref{eq:budget-law} describes the observed finite-budget range; it is not a complexity law. Together with the dimension comparisons above, it indicates sensitivity to training error without separating optimization and approximation errors.

Budget allocation changes which schedule is ahead. Below
$2.4\times10^4$ epochs, the $4000$-epoch schedule performs better:
at $B=1.6\times10^4$, its error is $3.88\times10^{-1}$ against
$4.55\times10^{-1}$. Above that budget, the $8000$-epoch schedule is
ahead: at $B=3.2\times10^4$, it gives $2.77\times10^{-1}$ against
$3.12\times10^{-1}$. The direct baseline does not show this improvement:
after $64000$ epochs, its error is $6.6\times10^{-1}$, within $2\%$
of its $16000$-epoch value.

\subsection{Experiment 3: training-schedule ablation}
\label{subsec:training-schedule-ablation}

The preceding comparisons use the same training budget at every Newton
step. We now fix the total budget and number of steps and vary the
allocation between earlier and later corrections. The residual indicator
\eqref{eq:numerical-certificate} involves only the final correction and
the final linearized residual: earlier training errors enter through the
current iterate.  We test this mechanism on the Cole--Hopf benchmark of
Section~\ref{subsec:exp-cole-hopf}, using single precision and no L-BFGS.
Three schedules allocate the same total budget of $8000$ Adam epochs over
$K=4$ Newton steps:
\[
\begin{aligned}
    \text{uniform}&=(2000,2000,2000,2000),\\
    \text{late}&=(500,500,1000,6000),\\
    \text{early}&=(6000,1000,500,500),
\end{aligned}
\]
each with three random seeds.

\begin{figure}[htbp]
    \centering
    \includegraphics[width=\textwidth]{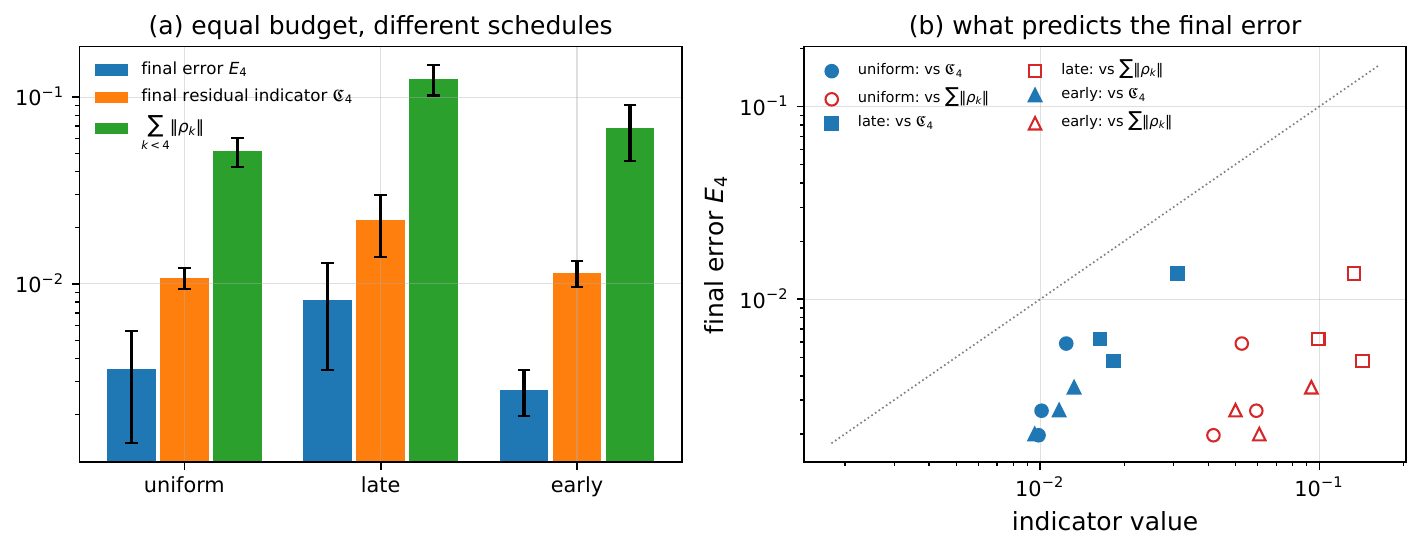}
    \caption{Experiment 3 (equal-budget schedule ablation, nine runs).
    (a) Final error, final residual indicator, and the sum of linearized residual norms
$\sum_{k=0}^{3}\norm{\rho_k}$ per schedule (mean $\pm$ standard deviation over
    seeds).  (b) Final error against the final indicator (filled) and
    against the accumulated residuals (open): only the final indicator
    tracks the final error in these runs.}
    \label{fig:ablation}
\end{figure}

Across all nine runs (Figure~\ref{fig:ablation}),
\[
    E_4/\mathfrak C_4\in[0.20,\,0.48],
\]
while the final errors span a factor of seven.  The final indicator tracks
the final error under every schedule; the accumulated quantity
$\sum_{k<4}\norm{\rho_k}$ 
exceeds $E_4$ by one to two orders of magnitude
and misorders the runs (Figure~\ref{fig:ablation}(b)).  

The same total budget does not give the same accuracy.  Mean final errors
are $2.7\times10^{-3}$ (early), $3.5\times10^{-3}$ (uniform), and
$8.2\times10^{-3}$ (late).  In these runs, a larger final-step budget did
not compensate for shorter training of the earlier corrections:
the mean $\norm{\rho_3}$ for the late schedule is about twice the corresponding means for the uniform and early schedules despite $6000$ final-step epochs, compared with $2000$ for the uniform schedule and $500$ for the early schedule.  

\subsection{Failure detection: growing interaction order}\label{subsec:exp-wide}

The next two benchmarks test whether the residual indicator remains
informative when training does not resolve the solution. First, a family with growing interaction order tests the prescribed architecture and training budget on increasingly broad phase interactions:
\begin{equation}\label{eq:wide-family}
    u^{\rm wide}_{\rm ex}(x,t)
    =
    \sin\Bigl(\sum_{i\in A}x_i\Bigr)
    +t\cos\Bigl(\sum_{i\in B}x_i\Bigr)
    +e^{-t}\sin\Bigl(\sum_{i\in C}x_i\Bigr),
\end{equation}
where $A$, $B$ are the first and last $\lceil n/2\rceil$ coordinates and $C$
the odd coordinates, so the phase interaction order is $\lceil n/2\rceil$. We use \eqref{eq:exp-hd-pde} with the Hamiltonian
\eqref{eq:exp-hd-H}, manufacturing the source from \eqref{eq:wide-family}. The protocol of the main sweep is run for $n\in\{5,10,20\}$ (three seeds).

\begin{figure}[htbp]
    \centering
    \maybegraphics{width=0.92\textwidth}{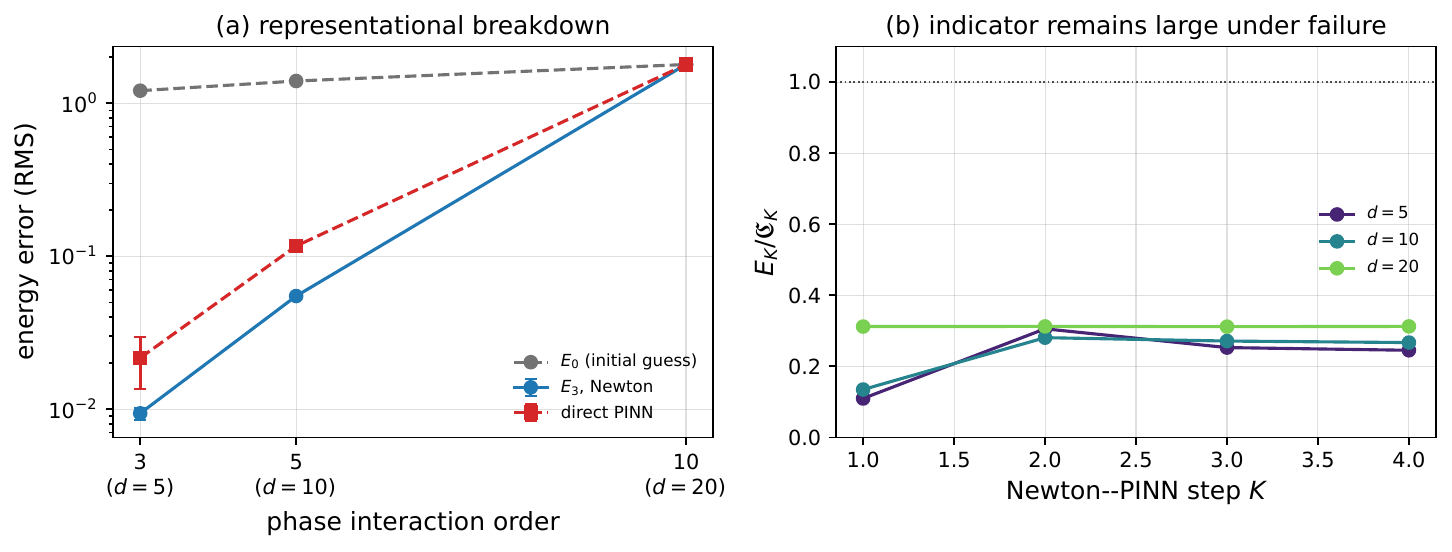}
    \caption{Wide family \eqref{eq:wide-family}: growing interaction order.
    (a)Errors at $K=0$ and $K=3$ (Newton, $12000$ epochs), and for the
direct PINN after $16000$ epochs.
    (b) Empirical calibration ratio along the iteration.}
    \label{fig:wide}
\end{figure}

At $n=5$ and $n=10$, respectively,
$E_3=9.4\times10^{-3}$ and $5.5\times10^{-2}$, against
$2.2\times10^{-2}$ and $1.2\times10^{-1}$ for the direct baseline after
$16000$ epochs, compared with $12000$ epochs for three Newton corrections. At $n=20$, the interaction order is $10$ and neither method appreciably improves on the initial energy error $E_0\approx1.79$ over the full $16000$-epoch budget. The indicator also fails to decrease: $\mathfrak C_K\approx5.7$ throughout, with $E_K/\mathfrak C_K=0.31$.

\subsection{An unforced high-dimensional benchmark}
\label{subsec:exp-unforced}

Every problem in Experiment 2 is forced.  This benchmark tests the
residual-monitoring pipeline without a manufactured source, using the
quadratic viscous problem
\[
\begin{cases}
    u_t+\tfrac12|Du|^2=\Delta u,\\
    u(x,0)=-2\log\phi_0(x).
\end{cases}
\]
The Cole--Hopf transformation gives an exact solution on $\T^n$ for any
smooth positive $\phi_0$.  We choose the nonseparable initial function
\[
    \phi_0(x)
    =2+\frac{a}{n}\sum_{i=1}^n\cos(x_i+x_{i+1})
    +\frac{b}{n}\sum_{i=1}^n\cos(x_i+2x_{i+1}-x_{i+2})
\]
with cyclic indices and $a=b=0.9$.  Its solution is
\[
\begin{aligned}
u_{\rm ex}(x,t)=-2\log\Biggl[
&2+\frac{a}{n}e^{-2t}\sum_{i=1}^n\cos(x_i+x_{i+1})\\
&+\frac{b}{n}e^{-6t}\sum_{i=1}^n
\cos(x_i+2x_{i+1}-x_{i+2})
\Biggr],
\end{aligned}
\]
with $M_2=1$ exactly and $S=0$.
The code uses the quadratic Hamiltonian itself, for which $D_pH(p)=p$ is
not globally bounded.  The stability proof applies with a bound on the
mean-value drift $(D\widehat u^K+Du_{\rm ex})/2$.  We do not evaluate this
bound and therefore report only the empirical unscaled indicator.

With the stated coefficients, the $1/n$ normalization gives
$\phi_0\ge 2-a-b=0.2$ uniformly in $n$. Consequently, the time-dependent signal has RMS $O(n^{-1/2})$ and $|Du|^2=O(1/n)$, and the problem approaches the heat equation as $n$ grows. Unlike the forced families, this family also has a vanishing signal scale. It therefore tests unforced residual monitoring with a dimension-dependent signal amplitude. The protocol of Experiment 2 is run at $n\in\{10,50,100\}$ (two seeds); training, stopping, and residual monitoring never use $u_{\rm ex}$.

At $n=10$, the iteration reduces the error from $E_0=0.55$ to
$E_4=5.4\times10^{-2}$, ahead of the direct baseline at the same epoch
count ($6.6$--$7.3\times10^{-2}$).  The calibration ratio
$E_K/\mathfrak C_K=0.29$--$0.43$ overlaps the lower end of the calibration range in Experiment 2. At $n=50$ and $n=100$,  fixed-budget training reduces the error by only $10\%$ and $2\%$:
$E_4/E_0=0.90$ and $0.98$, respectively.  The indicator reports this
stagnation, staying within $7\%$ of the initial residual for all $K$, with calibration ratio $0.39$--$0.43$. 

\subsection{Indicator-driven stopping}
\label{subsec:exp-adaptive}

The preceding experiments use the residual indicator for reporting.  The
final experiment tests whether it can also control the iteration, and
whether the same residual-based inner stopping rule is effective across
dimensions.

Each correction trains in blocks of $500$ epochs, with $\rho_k$ evaluated on a validation sample independent of the training batches.  Training stops when the inexact-Newton forcing condition
$\norm{\rho_k}\le\eta\,\norm{\mathcal R[\widehat u^k]}$, with $\eta=0.5$, holds or a per-step budget of $8000$ epochs is exhausted.  
The trained correction $h$ is accepted if $\norm{\mathcal R[\widehat u^k+h]}<\norm{\mathcal R[\widehat u^k]}$. Otherwise, we test $0.5h$ and then $0.25h$, both relative to the original trained correction, and accept the first trial that decreases the nonlinear residual norm. If neither trial does so, the correction is rejected and the iteration terminates at the last accepted iterate.
For every accepted step, $h^k$ denotes the accepted correction. We evaluate $\rho_k$ and $\mathfrak C_{k+1}$ using this correction.
The forcing test is applied before damping. Its condition may remain
unmet when the epoch cap is reached and need not hold after damping;
the recomputed indicator remains applicable.

The outer iteration terminates when $\mathfrak C_{k+1}\le\tau$,
when $k\ge1$ and $\mathfrak C_{k+1} \ge 0.98\,\mathfrak C_k$
(floor), or after six accepted corrections.
True errors are computed for reporting only: $\tau$ is an indicator
threshold in the sense of Section~\ref{subsec:protocol}, not a certified error tolerance. The controller is run on the ring family at
$n\in\{10,20,50\}$, with three seeds and $\tau\in\{0.3,0.1,0.05\}$ at $n\le20$ and $\tau=0.3$ at $n=50$. Thresholds below the observed fixed-budget indicator floor at $n=50$ are excluded.

\begin{table}[htbp]
    \centering
    \small
    \setlength{\tabcolsep}{5pt}
    % Auto-generated by make_phase123_tables.py -- do not edit by hand.
\begin{tabular}{llcccccc}
\toprule
$n$ & $\tau$ & stop reason & $K$ & epochs & $E$ & $\mathfrak C$ & $E\le\tau$ \\
\midrule
10 & 0.3 & indicator & 2 & 2500--4000 & $4.10(40){\times}10^{-2}$ & $1.02(11){\times}10^{-1}$ & 3/3 \\
10 & 0.1 & indicator & 2--3 & 3500--10500 & $3.53(63){\times}10^{-2}$ & $8.87(125){\times}10^{-2}$ & 3/3 \\
10 & 0.05 & indicator\_floor & 5--6 & 28000--35500 & $2.17(7){\times}10^{-2}$ & $6.29(9){\times}10^{-2}$ & 3/3 \\
\midrule
20 & 0.3 & indicator & 2 & 8500 & $1.01(2){\times}10^{-1}$ & $2.27(4){\times}10^{-1}$ & 3/3 \\
20 & 0.1 & max\_steps & 6 & 40500 & $6.49(9){\times}10^{-2}$ & $1.57(1){\times}10^{-1}$ & 3/3 \\
20 & 0.05 & max\_steps & 6 & 40500 & $6.49(9){\times}10^{-2}$ & $1.57(1){\times}10^{-1}$ & 0/3 \\
\midrule
50 & 0.3 & indicator\_floor & 5--6 & 33000--41000 & $1.99(4){\times}10^{-1}$ & $4.32(9){\times}10^{-1}$ & 3/3 \\
\bottomrule
\end{tabular}
    \caption{Indicator-driven stopping (ring family, three seeds per
    row).  Stop reason, corrections used, training epochs (range over
    seeds), final error and empirical indicator, and the fraction of runs
    whose true error lies below the requested indicator threshold.  The
    fixed-schedule protocol uses $16000$ epochs.}
    \label{tab:adaptive}
\end{table}

In all $9$ runs that terminated with $\mathfrak C_K\le\tau$, the true error also satisfies $E_K\le\tau$.  This is empirical calibration for one benchmark family, not a proof of coverage. Runs that do not reach the threshold stop at an indicator floor or the limit of six outer corrections.  At $n=10$, $\tau=0.3$ is reached after $2500$--$4000$ epochs, a reduction by a factor of $4$--$6.4$ relative to the fixed epoch budget; $\tau=0.1$ uses $3500$--$10500$ epochs.

The tested controller gives larger errors than the fixed schedule
in some cases. At $n=20$ (for $\tau\le0.1$) and $n=50$, the controller's
mean errors are $1.5$--$1.8$ times those of the fixed $K=4$ schedule at comparable or larger total budget. This is consistent with the schedule sensitivity in Experiment~3: shorter training of early corrections can leave errors that the tested later budgets do not recover.  

\subsection{Discussion}
\label{subsec:numerics-discussion}

 Three corrections reach accuracy
comparable to the grid solver at $n=1$.  The forced sweep
(Figure~\ref{fig:sweep-main}) reduces the initial energy error by
$89\%$--$98\%$ at $n=10,20,50$ and gives partial resolution at $n=100$.
The staged nonlinear and Newton baselines attain comparable accuracy in
the tested cases, while frozen linearization reduces training time and
yields the residual decomposition underlying the indicator.

The empirical indicator stays calibrated within factors $1.8$--$2.7$ for
$K\ge2$ in the main sweep and $2$--$5$ over the schedule ablation, and it
remains large under failure and stagnation.  These observations support
monitoring the current residual throughout training.  The main forced
sweep indicates that three Newton corrections are typically sufficient
before saturation for $n\le50$, whereas at $n=100$ either additional
corrections or an increased per-step budget remains beneficial. 
The observed indicator plateaus motivate an outer stopping test.
The adaptive results do not support applying the tested inner stopping
rule with $\eta=0.5$ unchanged across dimensions.

All of these conclusions retain the limitations stated above:
manufactured forcing in the main high-dimensional sweep, the memory-adapted timing protocol, the shared architecture of the component networks in C1--C3, and the unscaled empirical use of $\widehat{\mathfrak C}_K$.

%%%%%%%%%%%%%%%%%%%%%%%%%%%%%%%%%%

%%%%%%%%%%%%%%%%%%%%%%%%%%%%%%%%%%

%%%%%%%%%%%%%%%%%%%%%%%%%%%%%%%%%%%

\section*{Code and data availability}
Code and data for Section~\ref{sec:numerics} are available at
\url{https://github.com/yeoneung/newton_policy}.

\section*{
AI Statement}
The authors conceived, developed, and theoretically established all core ideas, mathematical formulations, and proofs presented in this work. 

The authors used ChatGPT, Gemini, and Claude solely to debug and refine code for numerical experiments and to edit the manuscript for grammar and clarity.
% for grammatical editing, clarity improvements, and language polishing of the manuscript. 
% The authors 
% % reviewed and edited the final content and 
% take full responsibility for the integrity and accuracy of the paper.
\bibliographystyle{abbrv}  
\bibliography{bib, mfgv10_nn} 

\end{document}